\documentclass[12pt]{article}

\usepackage{amsfonts}
\usepackage{color}
\usepackage{amssymb}
\usepackage{amsmath}
\usepackage{indentfirst}

\def\sqr#1#2{{\vcenter{\vbox{\hrule height.#2pt
				\hbox{\vrule width.#2pt height#1pt \kern#1pt \vrule width.#2pt}
				\hrule height.#2pt}}}}
\def\signed #1{{\unskip\nobreak\hfil\penalty50
		\hskip2em\hbox{}\nobreak\hfil#1
		\parfillskip=0pt \finalhyphendemerits=0 \par}}
\def\endpf{\signed {$\sqr69$}}

\def\dbR{{\mathop{\rm l\negthinspace R}}}

\def\3n{\negthinspace \negthinspace \negthinspace }
\def\2n{\negthinspace \negthinspace }
\def\1n{\negthinspace }

\def\ds{\displaystyle}

\def\dbR{{\mathop{\rm l\negthinspace R}}}
\def\={\buildrel \triangle \over =}

\def\a{\alpha}
\def\b{\beta}

\def\e{\varepsilon}

\def\l{\lambda}

\def\n{\nabla}

\def\t{\times}
\def\f{\varphi}
\def\th{\theta}
\def\o{\omega}

\def\i{\infty}
\def\ns{\noalign{\ss} }

\def\R{{\bf R}}

\def\D{\Delta}

\def\L{\Lambda}
\def\Si{\Sigma}

\def\O{\Omega}

\def\cA{{\cal A}}
\def\cB{{\cal B}}

\def\cH{{\cal H}}

\def\cL{{\cal L}}

\def\cP{{\cal P}}

\def\cl{{\cal l}}
\def\no{\noindent}

\def\ss{\smallskip}
\def\ms{\medskip}

\def\q{\quad}
\def\qq{\qquad}
\def\hb{\hbox}

\def\pa{\partial}

\def\cd{\cdot}

\def\div{\hbox{\rm div$\,$}}

\def\cl{\overline}

\def\|{\Big |}
\def\({\Big (}
\def\){\Big )}
\def\[{\Big[}
\def\]{\Big]}
\def\be{\begin{equation*}}
\def\bel{\begin{equation}\label}
\def\ee{\end{equation}}
\def\bt{\begin{theorem}}
\def\bcd{\begin{condition}}
\def\ecd{\end{condition}}
\def\et{\end{theorem}}
\def\bc{\begin{corollary}}
\def\ec{\end{corollary}}
\def\bde{\begin{definition}}
\def\ede{\end{definition}}
\def\bl{\begin{lemma}}
\def\el{\end{lemma}}
\def\bp{\begin{proposition}}
\def\ep{\end{proposition}}
\def\br{\begin{remark}}
\def\er{\end{remark}}
\def\ba{\begin{array}}
\def\ea{\end{array}}
\def\ed{\end{document}}
\def\ns{\noalign{\ms}}
\def\ds{\displaystyle}

\def\square#1{\vbox{\hrule\hbox{\vrule height#1%
			\kern#1\vrule}\hrule}}
\def\rectangle#1#2{\vbox{\hrule\hbox{\vrule height#1%
			\kern#2\vrule}\hrule}}

\font\tenbb=msbm10 \font\sevenbb=msbm7 \font\fivebb=msbm5

\newfam\bbfam
\scriptscriptfont\bbfam=\fivebb \textfont\bbfam=\tenbb
\scriptfont\bbfam=\sevenbb

\def\oO{{\overline \O}}

\newtheorem{lemma}{Lemma}[section]
\newtheorem{remark}{Remark}[section]

\newtheorem{theorem}{Theorem}[section]
\newtheorem{corollary}{Corollary}[section]

\newtheorem{definition}{Definition}[section]
\newtheorem{proposition}{Proposition}[section]
\newtheorem{condition}{Condition}[section]

\makeatletter

\@addtoreset{equation}{section}
\makeatother
\begin{document}
\title{\bf A unified weighted inequality for fourth-order partial differential operators and applications }

\author{Yan Cui\thanks{Department of Mathematics, Jinan University, Guangzhou 510632, China. cuiy32@jnu.edu.cn}, \  Xiaoyu Fu\ \thanks{School
		of Mathematics, Sichuan University, Chengdu
		610064, China.  xiaoyufu@scu.edu.cn }\  and \  Jiaxin Tian\thanks{School of Mathematics, Sichuan University, Chengdu 610064, China. tjx20102495@163.com }}

\date{}

\maketitle

\begin{abstract}
In this paper, we establish a fundamental inequality for  fourth order partial differential operator $\cP\=\a\pa_s+\b\pa_{ss}+\D^2$ ($\a, \b\in\dbR$) with an abstract exponential-type weight function. Such kind of weight functions including not only the regular weight functions but also the singular weight functions. Using this inequality we are able to prove some Carleman estimates for the operator $\cP$ with some suitable boundary conditions in the case of $\b<0$ or $\a\neq 0, \b=0$.
As application,  we obtain a resolvent estimate for $\cP$, which can imply a log-type stabilization result for the plate equation with clamped boundary conditions or hinged boundary conditions.
\end{abstract}

\no{\bf 2010 Mathematics Subject  Classification.}
 35G15, 35K35,93B05

\no{\bf Key Words}.  fourth order partial differential operators, plate equation, Carleman estimate.

\section{Introduction and main results}


 In 1939,  Carleman \cite{c}  introduced a method to prove the strong unique continuation property for second order elliptic partial differential equations in dimensional two. His method, now known as ``Carleman estimate". In 1963, H\"ormander \cite{Hormander63} provided  strong pseudo-convexity condition,  generalized this type of estimates to a large class of differential operators. Now, this method has become one of the major tools in the study of unique continuation property (e.g., \cite{Hormander, Zuily83}), inverse problems (e.g., \cite{Isakov2, Kn, Kli,  Y}) and control of partial differential equations (see \cite{FI, FLZ} and the rich references therein).

Compared to  second order partial differential operator, to derive Carleman-type inequality for higher order partial differential operator is not an easy task, especially for multi-dimensional case.  On the one hand, a classical method to obtain Carleman estimates for even order differential operator is that one can check couple $(\cP,\ell)$ satisfying  subelliptic conditions (see \cite{BL15} for details), where $\cP$ is the principal operator, $\ell$ is a smooth weight function. However, this is not always working for  differential operator of order not less than 4 and  a loss of one full derivative may occur, we refer to \cite{JL20} for details.
On the other hand,
boundary behavior of $\cP$ will be a big obstruction for the finer
investigation of Carleman estimates.

In the case of  higher order partial differential operators in one space dimension,
 Carleman estimates for Kuramoto-Sivashinsky equation (see \cite{CM, GP}), fourth order
parabolic equation  (see \cite{GP0, XCY, ZZ}) and fourth order Schr\"odinger operator (see \cite{Zheng}) were established. In the case of multi-dimensional case, there are only few references in this respect. In 2020, Le Rousseau and Robbiano \cite{JL20} studied spectral inequality and resolvent estimate for bi-Laplace operator $\D_g^2$ on a compact manifold. By virtue of micro-local analysis approach, the authors in \cite{JL20} established a local Carleman estimate for the bi-Laplacian operator with  clamped  boundary conditions. In 2019, Guerrero and Kassab \cite{GK} first established a global Carleman estimate for the fourth order parabolic-type operator $\pa_t+\D^2$ in dimension $N\ge 2$ with homogeneous Dirichlet boundary conditions on the solution and the Laplacian of the solution. Recently, Kassab \cite{Kassab} established a new Carleman estimate for the operator $\pa_t+\D^2$ with Dirichilet and Neumann boundary conditions, via which the author obtained null controllability of semi-linear fourth order parabolic equations. In 2022, Huang and Kawamoto generalized the estimate in \cite{GK}  to the operator $\pa_t-A^2$ with $A=\n\cdot (a(x)\n )$ and then obtained a stability result for fractional operator $\pa_t^{\frac{1}{2}}-A$ ( see \cite{HK}).

  In this paper, we aim at establishing a unified Carleman estimate for the operator $$\cP\=\a \pa_s+\b \pa_{ss}+\D^2$$ with exponential-type weight function $\th(s, x)=e^{\ell(s, x)}$, where $\alpha,\beta\in \R$ and $\ell$ is a smooth function which will be given later. As we know, Carleman estimates can be regarded as weighted energy estimate with exponential-type weight function. Elementary calculus can be enough to grasp the main idea of Carleman estimate (See \cite[Chapter 1]{FLZ}). To obtain the core, there are two important ingredients need be considered, i.e.,
\begin{itemize}
    \item the choice of weight function $\th$;
    \item the decomposition of principal operator $\th (\cP \th^{-1})$.
  \end{itemize}

 Let  $I:=(b_1,b_2)\subset \dbR$ be a bounded interval and $\O\subset \dbR^n$ be a bounded domain with a $C^4$ boundary $\pa\O$.  Set
 \bel{0329-0a} Q_{b_1,b_2}=(b_1, b_2)\t\O,\q \Sigma_{b_1,b_2}=(b_1,b_2)\times \pa\O.
 \ee

About the choice of weight function, we have the following assumption.
 \bcd\label{cd1}
There exists a function  $\xi\in C^{2,4}(Q_{b_1,b_2}; \dbR)$ satisfying the
following:
 \begin{enumerate}
 \item [{\rm(i)}] For any $\l>1$,  put
\bel{0417-0b}\ba{ll}\ds
\th=\th(s,x)= e^{\ell(s,x)},\q \ell=\l\xi=\l\xi(s,x).
 \ea\ee

 \item [{\rm(ii)}] For any $\mu>1$,  assume that there is a $\varphi=\varphi(s,x)>0$ and a function $\eta=\eta(x)$  such that
 \bel{0421-a}\left\{\ba{ll}\ds
\xi_{x_j}=\mu  \varphi \eta_{x_j},\q j=1,2,\cdots, n,\\
 \ns\ds \xi_{x_jx_k}=\mu^2\varphi\eta_{x_j}\eta_{x_k}+\mu\varphi\eta_{x_jx_k},\q j,k=1,2,\cdots, n
 \ea\right.\ee
 and there is a constant $C>0$ such that
  \bel{1211-0a0}
  |\varphi_s|\le C\varphi^3,\q |\varphi_{ss}|\le C\varphi^5.
  \ee
 \end{enumerate}
 \ecd
 This assumption actually is as usual in the context of $L^2$ Carleman estimates for second order operator.  There exist many functions satisfying Condition \ref{cd1}, for example,
 {\small
 \begin{tabular}{|c|c|c|c|p{4.2cm}|}
\hline
 & $I$ & $\xi$ &$\varphi$& $\eta$ \\
\hline
i)& $ b_1,b_2\in \dbR$ & $\xi=\frac{e^{\mu \eta}-2e^{\mu |\eta|_{C(\oO)}}}{{(s-b_1)}^k(b_2-s)^k}$ &$\varphi=\frac{e^{\mu \eta}}{{(s-b_1)}^k(b_2-s)^k}  (k\ge  \frac{1}{2})$& $\eta\in C^4(\oO) $\\

\hline
ii)& $I\subset \dbR$ & $\xi=e^{\mu\(\eta(x)-c(s-s_0)^2\)}$ &$\varphi=\xi$& $\eta=|x-x_0|^2, x_0\in \dbR^n\setminus\O$\\
\hline
iii)& $b=-b_1=b_2$ & $\xi=e^{\mu\({\eta(x)\over ||\eta||_{L^\infty(\O)}}+b^2-s^2\)}$ &$\varphi=\xi$&$\eta\in C^4(\oO) $\\
\hline

\end{tabular}
}

 It should be pointed out
 \begin{itemize}
   \item  weight functions presented in i)  are very useful
 for the controllability/observability problems of parabolic-type equation (see \cite{FI, GK} for example);
   \item weight functions presented in ii) can be used to studying the  controllability/observability problem for the plate equations or inverse problem (see \cite{Y} for example);
   \item weight functions presented in iii) are very useful for the logarithmic stabilization problems of partial differential equations (see \cite{FLZ, JL20} for example).
 \end{itemize}

 Now assume that $\th$ satisfies Condition \ref{cd1}. Then setting
 \bel{qwa}\ba{ll}\ds
 z=\th v.
 \ea\ee
Define
    \bel{2pde1}
        \ba{ll}
        \ds \cP v\=\a v_s+\b v_{ss}+\D^2 v,\q\a,\b\in\R.
        \ea
    \ee
 Let us consider the decomposition of operator $\th \cP (\th^{-1}z)$.
Notice that following from (\ref{0417-0b}) and (\ref{0421-a}), we have
\bel{0421-34}\left\{\ba{ll}\ds
 \ell_s=\l \xi_s ,\q \ell_{ss}=\l\xi_{ss}, \q \ell_{x_j}=\l\mu\varphi\eta_{x_j},\q j=1,2,\cdots, n,\\
 \ns\ds \ell_{x_jx_k}=\l\mu^2\varphi\eta_{x_j}\eta_{x_k}+\l\mu\varphi\eta_{x_jx_k},\q j,k=1,2,\cdots, n.
 \ea\right.\ee
 Hence, we get some derivative relations between $v$ and $z$
  \bel{0421-b}\left\{\ba{ll}\ds
 \th v_s= z_s-\ell_s  z,\q \th\n v=\n z-\n\ell v,\\
 \ns\ds \th v_{ss}= z_{ss}-2\ell_s z_s+(\ell_s^2-\ell_{ss})  z,\\
 \ns\ds \th\D v=\D z-2\n\ell\cdot\n z+(|\n\ell|^2-\D\ell)  z.
  \ea\right.\ee
For simplicity, we put
 \bel{0421-02}
A=|\n\ell|^2,\q \L=|\n\ell|^2-\D\ell
 \ee
 and
 $$
 \n^2 a:\n^2 b=\sum_{j,k=1}^n a_{x_jx_k} b_{x_jx_k},\q \n^2a\n c \n d=\sum_{j,k=1}^n a_{x_jx_k}c_{x_j} d_{x_k},
 $$
where $a,b,c, d\in C^2(\R^n;\R). $
Direct calculation yields
  \bel{0421-c}\ba{ll}\ds
  \th\D^2 v\\
  \ns\ds=\D^2 z-2\n\D\ell\cdot\n z-4\n^2\ell:\n^2 z-4\n\ell\cdot\n\D z+\D\L  z+2\n\L\cdot\n  z+2\L \D  z\\
  \ns\ds\q+4\n^2\ell\n\ell\n z+4\n^2 z\n\ell\n\ell-2\n\ell\cdot\n\L  z-4\L\n\ell\cdot\n z+\L^2 z.
  \ea\ee

Recalling $\cP$ in \eqref{2pde1} and the fact that $v=\th^{-1}z$, we have
\bel{as00}\ba{ll}\ds
\th\cP (\th^{-1}z)=\th\cP v=\th\(\a v_s+\b v_{ss}+\D^2 v\).
\ea\ee
Combining (\ref{0421-b}) and (\ref{0421-c}), we can divide  $\th\cP (\th^{-1}z)$ into
\bel{as0}\ba{ll}\ds
\th\cP (\th^{-1}z)=\cP_1(z)+\cP_2(z)+P_r(z),
\ea\ee
where
\bel{0227-01}\left\{\ba{ll}\ds
\vec{E}=4\n A+4\D\ell\n\ell,\q F=\Delta^2 \ell-2\n\cdot(A\n\ell),\\
\ns\ds
 \cP_1 (z)=\b z_{ss}+\D^2 z+A^2 z+2A\D z+4\n^2 z\n\ell\n\ell+\vec{E}\cdot\n z,\\
\ns\ds \cP_2 (z)=\a z_s-4 \n\ell\cdot\n\D z-4\n^2\ell:\n^2 z-2\D\ell\D z-4 A\n\ell\cdot\n z+ F z,
 \ea\right.\ee
  and
   \bel{0227-0-0}\ba{ll}\ds
  \cP_r (z)=\(\D\L+2\n\ell\cdot\n\D\ell-\D^2\ell+(\D\ell)^2-\a \ell_s+\b\ell_s^2-\b\ell_{ss}\) z\\
  \ns\ds\qq\qq +(2\n\L-4\n A-2\n\D\ell-4\D\ell\n\ell)\cdot\n  z -2\b\ell_s z_s.
  \ea\ee
We give some explanations for this decomposition here.
\begin{itemize}
    \item Operator $\cP_1$ is the symmetric part of  $\th\cP(\th^{-1})$. For any $y\in C_0^\i(\dbR^{1+n})$, it is easy to check that
 $$
 \int_{\dbR^{1+n}} \cP_1(z) ydxds=\int_{\dbR^{1+n}} z\cP_1(y)\mathrm{dxds}.
 $$
    \item Operator $\cP_2$ is the anti-symmetric part of  $\th\cP(\th^{-1})$. For any $y\in C_0^\i(\dbR^{1+n})$, it is easy to check that
 $$
 \int_{\dbR^{1+n}} \cP_2(z) ydxds=-\int_{\dbR^{1+n}} z\cP_2(y)\mathrm{dxds}.
 $$
  \item Operator $\cP_r$ is the rest lower order part of  $\th\cP(\th^{-1})$.
  \end{itemize}

Similar to \cite[Theorem 1.1]{FLZ}, multiplying $\th\cP(\th^{-1}z) $ by $\cP_1(z)$ can obtain that
 \bel{0227-f0}
 \th\cP(\th^{-1}z)\cP_1(z)=|\cP_1(z)|^2+\cP_2(z)\cP_1(z)+ \cP_r(z) \cP_1(z).
 \ee
 Noting that the left hand side (L.H.S. for short) of above equality is less than $|\th \cP v|^2+\frac{|\cP_1(z)|^2}{4}$, and $|\cP_1(z)|^2+\cP_r(z) \cP_1(z)$ is larger than $\frac{|\cP_1(z)|^2}{2}-\frac{|\cP_r(z)|^2}{2}$.
 So in order to establishing a weighted inequality for the fourth partial differential operators, the key step is to give the estimation of
 $\cP_2(z)\cP_1(z)$. On the one hand, the product of symmetric operator $\cP_1$ and the anti-symmetric operator $\cP_2$ will yield many good energy terms. On the other hand, there are some ``bad" terms arising. To modify the coefficients of the weighted energy terms such that they have a good sign, we introduced a lower order auxiliary operator $\cP_3$ to reduce the influence of these ``bad" terms. Suitable choice of $\cP_3$  will give us more flexibility in the proof of our main results. More precisely, we define
  \bel{01as}
 \cP_3 (z)= \Phi z+Z\n^2\eta: \n^2 z,
\ee
where  $\Phi$ and $Z$  are smooth functions (which will be given explicitly later).

Now, we rewrite (\ref{0227-f0}) as
 \bel{0227-0as0}\ba{ll}\ds
\th\cP(\th^{-1}z)\cP_1(z)=|\cP_1(z)|^2+\cP_2(z)\cP_1(z)+\cP_3(z)\cP_1(z)+ \cP_4(z) \cP_1(z),
\ea\ee
where
   \bel{0427-0}\ba{ll}\ds
  \cP_4 (z)=\cP_r(z)-\Phi z-Z\n^2\eta:\n^2 z.
  \ea\ee

Our main result can be stated as following:
\bt\label{point-th1}  Let $\a,\ \b\in \dbR$, $\th$ be given by (\ref{0417-0b}) satisfying Condition \ref{cd1}. Put $z=\th v$. Let $\Phi\in C^2(Q_{b_1,b_2};\dbR)$ and $Z\in C^2(Q_{b_1,b_2};\dbR^{n\times n})$ satisfying
 \bel{E}
 \Phi=-\b\lambda^3\mu^{\frac{7}{2}}\varphi^3|\n\eta|^4,\q Z= 8\lambda\mu\varphi,
  \ee
  where $\lambda,\mu,\varphi,\eta$ are given in Condition \ref{cd1}.
Then there exists a positive constant $\mu_0$,  such that for any $\mu \ge \mu_0$, one can find three positive constants $\l_0=\l_0(\mu)$, $c>0$ and $C = C(\mu)$,  such that for all $v\in C^{2,4}(Q_{b_1,b_2})$ satisfying
  $\cP v=\(\a\pa_s +\b\pa_{ss}+\D^2\) v(s,x)$, for all $\l\ge \l_0$, it holds that
  \bel{point-th1-01}\ba{ll}\ds
 c\th^2\(\lambda^6\mu^8\varphi^6|\n\eta |^8 v^2+\lambda^4\mu^6\varphi^4 |\n\eta |^6 |\n v |^2+\l^3\mu^4\varphi^3|\n\eta|^4\th^2|\D v |^2\)\\
  \ns\ds\q+c\lambda\mu^2\varphi|\n\eta\cdot\n\D z|^2+\b^2\lambda^3\mu^{\frac{7}{2}}\varphi^3|\n\eta|^4 z_s^2+8\b \l\mu\varphi\n^2\eta\n z_s\n z_s\\
  \ns\ds \q+64\lambda\mu\varphi\n^2\eta(\n^2 z\n\ell)(\n^2 z\n\ell)+3\lambda^5\mu^6\varphi^5|\n\eta|^4|\n\eta\cdot\n z|^2+  \pa_s M + \div  V\\
 \ns\ds\le C|\th \cP v|^2\\
 \ns\ds\q+C\th^2\[\l^4\mu^5\varphi^4|\n v|^2+\l^6\mu^7\varphi^6v^2+\l^2\mu^3\varphi^2|\n^2 v|^2+\b^2\l^2\mu^3\varphi^2v_s^2+\mu|\n\D v|^2\]
 \ea\ee
 where  $M$ and $V$ are given by (\ref{mv-00a0})  and \eqref{boundary-0a}, respectively.
 \et

\br
Since we did not give boundary conditions in Theorem \ref{point-th1}, the forms of $M$ and $V$ would be quite ``complex", but once appropriate boundary conditions are imposed, many terms in $M$ and $V$ would vanish.
\er

 Theorem \ref{point-th1} gives a weighted  estimate for $v$ and derivative of $v$ in terms of $\cP v$ and low order terms (the last term in the R. H. S of (\ref{point-th1-01}). Specially, if $|\n\eta|>0$ in $\oO$, these lower order terms can be absorbed by the L. H. S of (\ref{point-th1-01}). To prove Theorem \ref{point-th1}, the key point is the estimations of $\cP_1(z)\cP_2(z)$ and $\cP_1(z)\cP_3(z)$. We give some explanations here.
  \begin{itemize}
    \item {\bf Estimation of $\cP_1(z)\cP_2(z)$. } Using divergence formula, we can derive a pointwise estimate on the product of   $\cP_1(z)$  and $\cP_2(z)$ in terms of the sum of ``divergence" terms, ``energy" terms and some undesired lower order terms (see Lemma \ref{lem:1} and Lemma \ref{lem:02} in Section \ref{section1}).

    \item {\bf Estimation of $\cP_1(z)\cP_3(z)$. } To deal with some undesired terms appeared in the estimation of $\cP_1(z)\cP_2(z)$, we choose  suitable auxiliary functions $\Phi$ and $Z$ in $\cP_3(z)$, such that some positive terms would appear, which can help us to control those undesired terms (see Proposition \ref{lem:4} in Section \ref{section2} and Lemma \ref{lem:8} in Section \ref{section2}). As we see in (\ref{E}), $\Phi=0$ if $\b=0$. While in the case of $\b\neq0$, $\Phi$ will play an important role in our proof. Moreover, the choice of $\Phi$ can be relaxed to $\Phi=-\b\l^3\mu^\kappa\varphi^3|\n\eta|^4$ for any $\kappa\in(0,4)$. For simplicity, here we choose $\Phi=-\b\l^3\mu^{7/2}\varphi^3|\n\eta|^4$. It should be pointed out that the choice of  $\Phi$ can improve the coefficients of $z_s^2$. In this situation,   $-2\b\ell_s z_s$ changed to be a lower order term, hence we put it in the rest part in (\ref{0227-0-0}). Moreover, it does not influence the sign of $\cB_j (j=1, 2, 3)$ in (\ref{1209-f1}). The choice of $Z$ plays a key role in the estimation of $H_j (j=1, 2, 3)$ in (\ref{12101-a}).
  \end{itemize}

In the sequel, we give some applications of Theorem \ref{point-th1}. We consider $\cP$ with clamped boundary conditions \bel{bc:cbc}\ds
v=\frac{\pa v}{\pa\nu}=0\q \mbox{ on } \Si_{b_1,b_2},
  \ee
 or hinged boundary conditions,
  \bel{bc:hbc}\ds
   v=\D v=0\q \mbox{ on } \Si_{b_1,b_2}.
  \ee

Let $\o_0, \o$ be non-empty open subsets of $\O$ satisfying ${\cl{\o_0}}\subset\o$. According to \cite{FI} or \cite{coron}, we know that there exists a function   $\eta\in C^4(\oO)$ satisfying
 \bel{1224-0a}\ba{ll}\ds
\eta> 0\q \mbox{ in }\O, \qq\eta=0\q\mbox{ on }\pa\O,\qq|\n\eta|>0\q \forall x\in \overline{\O\backslash \omega_0}.
 \ea\ee

  We have the following Carleman estimate.
\bt\label{thm:2}
Let $\b<0$, $\th$  satisfying Condition \ref{cd1} and $\eta$ satisfying (\ref{1224-0a}). Then there exists a positive constant $\mu_0$,  such that for any $\mu \ge \mu_0$, one can find three positive constants $\l_0=\l_0(\mu)$, $c_1>0$ and $C = C(\mu)$,  such that for all $v\in H^{2, 4}(Q_{b_1,b_2})$ satisfying
 \bel{0305-a}
 v(b_1,\cdot)=v(b_2,\cdot)=v_s(b_1,\cd)=v_s(b_2,\cd)=0,
 \ee
 with clamped boundary conditions \eqref{bc:cbc} or hinged boundary conditions \eqref{bc:hbc}, for all $\l\ge \l_0$, it holds that
  \bel{i1i2-c1}\ba{ll}\ds
 \int_{Q_{b_1,b_2}}\th^2 \(\lambda^6\mu^8\varphi^6  v^2+\lambda^4\mu^6\varphi^4  |\n v |^2+\l^3\mu^\frac{7}{2}\varphi^3 (|\D v |^2+|v_s^2|\)\mathrm{dxds}
\\
  \ns\ds\q
+\int_{Q_{b_1,b_2}}\th^2 \l\mu^{\frac{3}{2}}\varphi\(|\n\D v|^2+|\n v_s|^2\)\mathrm{dxds}\\
 \ns\ds\le C\[||\th \cP v||^2_{L^2(Q_{b_1,b_2})}+\lambda^7\mu^{8}\int_{b_1}^{b_2}\int_\o\varphi^7\th^2 v^2\mathrm{dxds}\].
 \ea\ee
\et

\br
In the case of $\b<0$, Theorem \ref{thm:2} holds for $\a=\a(s,x)$.
\er

Based on  Theorem \ref{thm:2}, we can obtain an important resolvent estimate for the bi-Laplace operator. More precisely, we define an unbounded operator
 $$
 \cA\=\(\ba{ll}\ds 0, &1\\
 \ns\ds -\D^2, &-d \ea\)
 $$
 acting on $H_0^2(\O)\t L^2(\O)$ with domain $D(\cA)=(H^4(\O)\cap H_0^2(\O))\t H_0^2(\O)$, where $d=d(x)$ is a nonnegative function.
\bt\label{0tt2} Let $ \o$ be any fixed non-empty
open subsets of $\O$. Let $d_0, d_1$ be any
fixed positive constants. Assume that $d(\cd)$ satisfy
\begin{equation}\label{cond1.1}
 0\le d(x)\le d_1 \ \hbox{ in } \O, \q    d(x)\ge d_0>0 \ \hbox{ in }\o.
\end{equation}
 Then the unbounded operator $\cA-i\gamma  I$ is invertible on $\cH=H_0^2(\O)\t L^2(\O)$ for all $\gamma\in\dbR$ and there exists a constant $ C>0$ such that
  \bel{0324-0f}
 |\!|(\cA-i\gamma  I)^{-1}|\!|_{\cL(\cH)}\le Ce^{C|\b|},  \q \forall  |\gamma|>1.
  \ee
 \et
 \br
Based on Frequency domain approach \cite{BZ, burq98} and \cite[Proposition
7.2]{D} , Theorem \ref{0tt2} can imply a log-type decay result for all  solutions to the following damped plate equation
\bel{1215-0a1}\left\{\ba{ll}\ds
   z_{tt}+\D^2 z+d(x) z_t=0 &\hb{ in } \dbR^+\t \O,\\
\ns\ds
(z(0),z_t(0))=(z^0,z^1)&\hb{ in } \O,
 \ea\right.\ee
with clamped boundary conditions $z=\frac{\pa z}{\pa\nu}=0\q \mbox{ on } \R_+\times \O$ or hinged boundary conditions, $v=\D v=0\q \mbox{ on } \R_+\times \O.$
More precisely, let $d(\cd)$  satisfy (\ref{cond1.1}). Then, the
solution $e^{t\cA}(z^0,z^1) \in C(\dbR;
D(\cA))\cap C^1(\dbR; \cH)$ to the system
(\ref{1215-0a1}) satisfies
\bel{1215-0a3}\ba{ll}\ds
|\!|e^{t\cA}(z^0,z^1)|\!|_{\cH}\le \frac{C}{\log (2+t)}|\!|(z^0,z^1)|\!|_{D(\cA)},\qq\forall\;(z^0,z^1)\in
D(\cA),\ \forall\; t>0. \ea\ee
 \er

In the case of  $\b=0,\a\neq0$, we have  the following result:

\bt\label{thm:3}
Let  $\a\neq 0, \b=0$, $\th$ satisfying Condition \ref{cd1} and $\th(b_1,\cdot)=\th(b_2,\cdot)=0$. Assume that $\eta\in C^4(\oO)$ satisfying
 (\ref{1224-0a}). Then there exists a positive constant $\mu_0$,  such that for any $\mu \ge \mu_0$, one can find three positive constants $\l_0=\l_0(\mu)$, $c_1>0$ and $C = C(\mu)$,  such that for all $v\in H^{1, 4}(Q_{b_1,b_2})$  satisfying
  the boundary conditions (\ref{bc:cbc}) or (\ref{bc:hbc}), for all $\l\ge \l_0$, it holds that

  \bel{i1i2-c10}\ba{ll}\ds
\int_{Q_{b_1,b_2}} \frac{1}{\l \varphi}\th^2\(|\a v_s|^2+|\D^2 v|\)\mathrm{dxds}\\
  \ns\ds\q+
\int_{Q_{b_1,b_2}}\th^2 \(\lambda^6\mu^8\varphi^6 v^2+\lambda^4\mu^6\varphi^4  |\n v |^2+\l^3\mu^4\varphi^3 |\D v |^2+\lambda\mu^2\varphi|\n\D v|^2\)\mathrm{dxds}\\
 \ns\ds\le C\[\|\|\th (\a v_s+\D^2 v)\|\|^2_{L^2(Q_{b_1,b_2})}+\lambda^7\mu^{8}\int_{b_1}^{b_2}\int_\o\varphi^7\th^2 v^2\mathrm{dxds}\].
 \ea\ee
\et
\br Similar to \cite{GK}, based on Hilbert Uniqueness Method (see \cite{Lionscontrolbook} or \cite{coron}), Theorem \ref{thm:3} can imply a null controllability result for the following control system:
\bel{system:control}\left\{\ba{ll}\ds
   y_{t}+\D^2 y=\chi_\omega u &\hb{ in } \dbR^+\t \O,\\
   \ns\ds  y=\frac{\pa y}{\pa \nu}=0   &\hb{ on } \dbR^+\t \pa\O,\\
\ns\ds
y(0,\cdot)=y^0(\cdot) &\hb{ in } \O,
 \ea\right.\ee
\er
where $u\in L^2(\Omega)$ is the control, $y^0\in L^2(\O)$ is the initial state.

The rest of this paper is organized as follows. We give some preliminaries we needed in Section \ref{section1}.  The proof of Theorem \ref{point-th1} is given in Section \ref{section2}. In Section \ref{section4-0} and Section \ref{section5-0}, we give the proofs of Theorem \ref{thm:2} and Theorem \ref{thm:3}, respectively. Finally, in Section \ref{section6-0}, we give the proof of resolvent estimate, that is Theorem \ref{0tt2}.


\section{Preliminaries}\label{section1}

In this section, we collect some preliminaries we needed. Throughout this paper, we always assume that $C$ denotes a positive constant that is independent with $\l$ and may vary from one place to another.

First, we define the low-order terms  in the right hand of \eqref{point-th1-01} as:
   \bel{A11-12-13}\left\{\ba{ll}\ds
A_1\=\lambda^3\mu^5\varphi^3(\mu+\lambda\varphi)\(|\n z|^2+\l^2\mu^2\varphi^2z^2+\lambda\varphi|\n\eta\cdot\n z|^2\), \\
\ns\ds
A_2\=\lambda\mu^3\varphi\(\mu+\lambda\varphi\)|\n^2 z|^2,\q  A_3=\b^2\l^2\mu^3\varphi^2|z_s|^2,\\
\ns\ds A_4=\mu|\n\D z|^2+\l\mu\varphi|\n\eta\cdot\n\D z|^2.
 \ea\right.\ee
 Recalling that $\ell$ satisfies Condition \ref{cd1} and $z=\th v$, it is easy to check that the low-order terms in \eqref{point-th1-01} holds
  \bel{Aii}\ba{ll}
\l^4\mu^5\varphi^4|\n v|^2+\l^6\mu^7\varphi^6v^2+\l^2\mu^3\varphi^2|\n^2 v|^2+\b^2\l^2\mu^3\varphi^2v_s^2+\mu|\n\D v|^2\ds\\
 \ns\ds\q \ge A_1+A_2+A_3+A_4.
	 \ea\ee
 Next, we give the estimations of  $\cP_4(z)\cP_1(z)$.
\bl \label{lem:Ai} There exists a constant $C>0$, such that
\bel{0301-as0}
    \cP_4(z)\cP_1(z)\geq -\frac{|\cP_1(z)|^2}{4} -C(A_1+A_2+A_3).
\ee
\el
{\bf Proof of Lemma \ref{lem:Ai}. } Recalling (\ref{0427-0}) for the definition of $\cP_4(z)$, by (\ref{0227-0-0}), (\ref{E}) and (\ref{0421-34})--(\ref{0421-02}), it is easy to see that
\bel{0427-22}\ba{ll}\ds
  \cP_4(z)\le C\lambda^2\mu^3\varphi^2(\mu+\lambda\varphi)| z|+C\lambda\mu^3\varphi |\n\eta\cdot\n z|+C\lambda\mu\varphi|\n^2 z:\n^2\eta|+C\ell_s| z_s|
  \ea\ee
  Thus by using inequality $ab\geq -\frac{a^2}{4}-b^2 $ and \eqref{A11-12-13}, we get the desired result immediately.
  \endpf
Hence, we only need consider the estimation $\ds\cP_1(z)\cP_2(z)+\cP_1(z)\cP_3(z)$.  For simplicity, we set
 \bel{0624-01a}\left\{\ba{ll}\ds
 J_1(z)=\D^2 z+A^2 z+2A\D z+4\n^2 z\n\ell\n\ell+\vec{E}\cdot\n z,\\
 \ns\ds J_2(z)=-4 \n\ell\cdot\n\D z-4\n^2\ell:\n^2 z-2\D\ell\D z-4 A\n\ell\cdot\n z+Fz.
 \ea\right.\ee
 Recalling (\ref{0227-01}) the definitions of $\cP_j(z)$ ($j=1,2$) and (\ref{01as})  the definition of $\cP_3(z)$, it is easy to see that
   \bel{1208-cp}\ba{ll}\ds
   \(\th\cP v-\cP_4(z)\)\cP_1(z)\\
   \ns\ds=|\cP_1(z)|^2+\a\b z_sz_{ss}+\a z_sJ_1(z)+\b z_{ss}J_2(z)+\b z_{ss}\cP_3(z)+J_1(z)J_2(z)+J_1(z)\cP_3(z).
   \ea\ee
As we explained before, the estimations of the right hand side of (\ref{1208-cp})  will yields many boundary terms. Here, all the boundaries can be divided into three parts.

\ss

\ss

{\bf Part I: Boundaries terms of ``$\ds\a\b z_sz_{ss}+\a z_sJ_1(z)+\b z_{ss}J_2(z)+\b z_{ss}\cP_3(z)$".}

   \bel{mv-00a0}\left\{\ba{ll}\ds
   M=M_1+M_2+M_3,\\
   \ns\ds
  M_1=\frac{1}{2}\(\a\b z_s^2+\a|\D z|^2+\a A ^2 z^2\)-A\a |\n z|^2-2\a |\n\ell\cdot\n z|^2,\\
  \ns\ds M_2=4\b\n\ell\cdot\n z_s\D z+4\b\D\ell\D z z_s+4\b \n\D\ell\cdot\n z z_s-\b \D\ell_s|\n z|^2-2\b z_s\D\ell\D z\\
 \ns\ds\qq\q-2\b \n^2\ell_s\n z\n z+4\b\n^2\ell\n z\n z_s-4\b  A\n\ell\cdot\n z z_s+\b F z_s z-\frac{1}{2}\b F_s z^2 ,\\
  \ns\ds M_3=\b\Phi z_s z-\frac{1}{2}\b \Phi_s z^2+\b Z z_{s}\n^2 \eta:\n^2 z+\frac{1}{2}\b Z_s\n^2 \eta\n z\n z,
  \ea\right.\ee
  and
   \bel{m0-a0}\left\{\ba{ll}\ds
  V_1=\a z_s\n\D z-\a\n z_s\D z+2\a A z_s\n z+4\a z_s(\n\ell\cdot\n z)\n\ell,\\
  \ns\ds V_2=-4\b \D z (z_{ss} \n\ell+\n\ell_sz_s)-4\b\(\n\ell\cdot\n z_s\)\n z_s+2\b\n\ell|\n z_s|^2-4\b\D\ell\n z_s z_s\\
  \ns\ds\qq\q-4\b\n^2\ell\n z z_{ss}+2\b z_s\D\ell\n z_s+2\b \D\ell_s\n z z_s+\b (2  A\n\ell-\n\D\ell) z_s^2,\\
  \ns\ds V_3=-\b  Z_s z_{s}\n^2 \eta\n z-\b Z z_{s}\n^2 \eta\n z_s+\frac{1}{2}\b (\n^2 \eta \n Z) z_s^2+\frac{1}{2}\b Z \n\D \eta z_s^2.
   \ea\right.\ee
Obviously, in the case of $\b=0$, $M_2=M_3=V_2=V_3=0$. In the case of $\ds z=\frac{\pa z}{\pa\nu}=0$ on the boundary, it is easy to see that $V_1=V_2=V_3=0$.

   \ss

 {\bf Part II: Boundaries terms of ``$\ds J_1(z)J_2(z)+J_1(z)\cP_3(z)$".}
       \bel{m0v-a0}\left\{\ba{ll}\ds
   V_{4}=-4\n\D z(\n\ell\cdot\n\D z)+2\n\ell|\n\D z|^2-4\n\D z(\n^2\ell:\n^2 z)+4\D z\cdot\n A(\n\ell\cdot\n z)\\
  \ns\ds\qq\q+4\n^2 \ell\cdot(\n^2 z\cdot\n\D z)-4\n^2 z\cdot(\n^2\ell\cdot\n\D z)-2\D\ell\D z\n\D z\\
  \ns\ds\qq\q-4\n\D z(A\n\ell\cdot\n z)+ 4A\n^2\ell\n z\D z+ 4A \n^2 z\n\ell\n^2 z-2A\n\ell |\n^2 z |^2\\
  \ns\ds\qq\q-4\D z\n z\cdot\n A\n\ell+4\n^2 z\n\ell\n z\cdot\n A-4A\n^2\ell(\n z\cdot\n^2 z)\\
\ns\ds\qq\q+2\n A\cdot\n^2\ell|\n z|^2+4A\n^2\ell:\n^2 z\n z+F \n\D z z-F\D z\n z,\\
\ns\ds\qq\q -4\n^2 z \n\D\ell\D z+4Az \n^2z \n\D\ell\\
\ns\ds V_{5}= -4 A^2 z\n^2 z\n\ell+2A^2\n\ell |\n z|^2-2A^2\D\ell z\n z-2A^3 \n\ell z^2+4 z\n A^2\(\n z\cdot \n\ell\)\\
\ns\ds\qq\q-4A\n\ell|\D z|^2-8A^2\n z(\n\ell\cdot\n z)+4A^2\n\ell |\n z|^2+2A F z\n z,\\
\ns\ds V_{6}=-16\(\n^2 z\n\ell\n\ell\)\n^2 z\n\ell+8| \n^2 z\n\ell |^2\n\ell-8A\n\ell|\n\ell\cdot\n z|^2-8\D\ell\D z(\n z\cdot\n\ell)\n\ell\\
 \ns\ds\qq\q+8\D\ell\(\n^2 z\cdot\n\ell\)\(\n z\cdot\n\ell\)-4(\vec{E}\cdot\n z)\n^2 z\n\ell+4F z(\n z\cdot\n\ell)\n\ell,\\
  \ns\ds
 V_{7}=\Phi \n\D z z-\Phi\D z\n z+( Z\n^2 z:\n^2\eta)\n\D z+A^2{ Z}(\n^2\eta\n z) z-Z\n^2\eta\( \n^2  z\cdot\n\D z\) \\
 \ns\ds\qq+Z\n^2 z\( \n^2 \eta\cdot\n\D z\)+2A\Phi z\n z+4\Phi z(\n z\cdot\n\ell)\n\ell.
  \ea\right.\ee

\ss

{\bf Part III: The rest boundaries arising in the process of retrieving the lost energy terms. }
  \bel{mv-0331a}\left\{\ba{ll}\ds
  V_{8}=4(\nabla\ell\cdot\nabla {u})\nabla u-2|\nabla u|^2\nabla\ell+2\D\ell{u}\nabla u+2|\nabla\ell|^2\nabla\ell|u|^2,\\
 \ns\ds
  V_{9}=32\lambda\mu\varphi\n\ell(\n^2 z:\n^2\eta)\n z\cdot\n\ell-32\lambda\mu\varphi\n^2\eta(\n z\cdot\n\ell)(\n^2 z\n\ell)\\
  \ns\ds\qq-16A\lambda\mu^2\varphi\n\eta\n z\cdot\n\eta\D z-8A\n A\cdot\n\ell z\n z+2A\b\lambda^3\mu^{\frac{7}{2}}\varphi^3|\n\eta|^4 z\n z,
   \ea\right.\ee
   where  $u=\varphi^{3/2}|\n\eta|^2v$.

 Define
  \bel{boundary-0a}\ba{ll}\ds
  V=\sum_{j=1}^8 V_j-\l^3\mu^4V_8+V_9.
  \ea\ee
 \bl\label{lem:1} Let $\a, \b\in\dbR$, $J_k (k=1, 2), \cP_3(z)$ be given by (\ref{0624-01a}) and (\ref{01as}), respectively. Assume that $\ell\in C^4(Q_{b_1,b_2})$ satisfying Condition \ref{cd1}. Let $\Phi\in C^2(Q_{b_1,b_2};\dbR)$ and $Z\in C^2(Q_{b_1,b_2};\dbR^{n\times n})$. 
For any $z\in C^{2,4}(\dbR^{1+n})$, we have
  \bel{i1i21-p1p2}\ba{ll}\ds
 \a\b z_sz_{ss}+\a z_sJ_1(z)+\b z_{ss}J_2(z)+\b z_{ss}\cP_3(z)\\
 \ns\ds
 \ge  \pa_s M+\div\(V_1+V_2+V_3\)+\b Z\n^2 \eta\n z_{s}\n z_s\\
 \ns\ds\q-\b\Phi z_s^2-C\(A_1+A_2+A_3+A_4\)
   \ea\ee
   \el
  where $M$ and $ V_j, (j=1, 2, 3)$ are given by \eqref{mv-00a0} and \eqref{m0-a0}, respectively.

Similarly, we have the following ``Estimation of all space-derivative  terms".
\bl\label{lem:02} Under the same conditions of Lemma \ref{lem:1},  we have
 \bel{1209-b}\ba{ll}&\ds
J_1(z)J_2(z)+J_1(z)\cP_3(z)\\
\ns&\ds\ge\div (V_4+V_5+V_6+V_7+V_8)+\cB_1|\n z|^2+\cB_2|\D z|^2+\cB_3 z^2+\cB_4|\n\ell\cdot\n z|^2\\
  \ns&\ds\q+(\cH_1+\cH_2+\cH_3+\cH_4)+\mathcal{R}- C\(A_1+A_2+A_3+A_4\)
\ea\ee
where $V_j (j=4, \cdots, 8)$ are given by \eqref{m0v-a0},
 \bel{1209-f1}\left\{\ba{ll}\ds
 \cB_1=-8A\n A\cdot\n\ell-2A\D^2\ell-2A\Phi,\\
 \ns\ds\cB_2=4\n A\cdot\n\ell+\D^2\ell+\Phi,\\
 \ns\ds\cB_3=4A^2\n A\cdot\n\ell+A^2\D^2\ell+A^2\Phi,\\
 \ns\ds\cB_4=8A\lambda\mu^2\varphi|\n\eta|^2+16\n A\cdot\n\ell-4\D^2\ell-4\Phi,
 \ea\right.\ee
and
 \bel{12101-a}\left\{\ba{ll}\ds
 \cH_1=8\n^2\ell\n\D z\n\D z-Z\n^2\eta\n\D z\n\D z,\\
 \ns\ds
  \cH_2=8A^2\lambda\mu\varphi\n^2 \eta\n z\n z-A^2{ Z}\n^2\eta\n z\n z,\\
 \ns\ds \cH_3=32\n^2\ell(\n^2 z\n\ell)(\n^2 z\n\ell)+4 Z(\n^2 z:\n^2\eta)(\n^2 z\n\ell\n\ell),\\
 \ns\ds\cH_4=-16A \n^2\ell:\n^2 z\D z+2A Z(\n^2 z:\n^2\eta)\D z.
 \ea\right.\ee
 And moreover, the remainder term $\mathcal{R}$ is given by
  \bel{1210-R4}\ba{ll}\ds
 \mathcal{R}=-\n \Phi \cdot\n\D z z+\n \Phi\cdot\n z \D z-2\n\(A\Phi\)\cdot \n z z +\Phi\vec{E}\cdot\n z z\\
 \ns\ds\qq\ -\(\n^2\eta\n(A^2{ Z})\)\n z z-A^2{ Z}(\n\D\eta\cdot\n z) z+ Z\(\n^2 z:\n^2\eta\)(\vec{E}\cdot\n z)\\
 \ns\ds\qq\ -4\(\n\Phi\cdot\n\ell\)(\n z\cdot\n\ell) z-4\Phi\D\ell(\n z\cdot\n\ell) z-4\Phi(\n^2\ell\n\ell\n z) z\\
\ns\ds\qq\ -\sum_{i, j, k=1}^n\[\(Z\eta_{x_ix_k}\)_{x_j}\D z_{x_j}z_{x_ix_k}-\( Z\eta_{x_ix_k}\)_{x_k}\D z_{x_j}z_{x_ix_j}+\( Z\eta_{x_ix_k}\)_{x_j}\D z_{x_k}z_{x_ix_j}\].
 \ea\ee
  \el

  \br
  To see the role of $\Phi$ and $Z$ clearly, we only assume $\Phi$ and $Z$ are smooth functions, here we haven't give the choice of them, so we have the remainder term $\mathcal{R}$. Further, if we choose $\Phi$ and $Z$ satisfying (\ref{E}), then it is easy to check that $\ds |R|\le C(A_1+A_2+A_4)$ (see (\ref{R1}) for details).
  \er

The proofs of Lemma \ref{lem:1} and Lemma \ref{lem:02} are elementary and long, so we out them in the Appendix.

To retrieve the lost energy terms, we recall the following known result.
 \bl\label{2l1} { \upshape(see \cite{FLZ})} Assume that $v\in C^2(\dbR^{1+n}; \;\dbR)$. Let $\th$ be given by (\ref{0417-0b}). Put
   \bel{0303-a01}
   u=\varphi^{3/2}|\n\eta|^2v.
   \ee
   Then there is a constant $C>0$ such that
 \bel{2a2}\ba{ll}
2\th^2|\D v|^2+\n\cdot V_8&\ds\ge
 2\lambda^3\mu^4\varphi^3|\nabla\eta|^4|v|^2 +2\lambda\mu^2\varphi|\nabla\eta|^2|\n v|^2\\
 \ns&\ds\q-C\[(\lambda^3\mu^3\varphi^3+\l^2\mu^4\varphi^2)|v|^2+\l\mu\varphi|\n v|^2\],
	 \ea\ee
where $ V_8$ is given by \eqref{mv-0331a}.
 \el

 Finally, we have the following result.
\bl\label{lem:5} Let  $\b<0$ and $\eta\in C^4(\oO)$ satisfying (\ref{1224-0a}). Fix $b_1, b_2\in\dbR$ with $b_1<b_2$. For all $v\in H^{2, 4}([b_1,b_2]\t\O)$ satisfying $v_s(b_1,\cd)=v_s(b_2,\cd)=0$.

i) In the case $v=\frac{\pa v}{\pa\nu}=0$ on the boundary, we have
\bel{2p14-0}\ba{ll}\ds
\int_{Q_{b_1,b_2}}\(-\b\lambda\mu^{\frac{3}{2}}\varphi \th^2|\n\D v|^2+\b^2 \lambda\mu^{\frac{3}{2}}\varphi \th^2 |\n v_s|^2\)\mathrm{dx} \mathrm{ds} \\
\ns\ds\leq C\[||\th\cP v ||^2_{L^2(Q_{b_1, b_2})}+\int_{Q_{b_1,b_2}}\lambda^3\mu^{\frac{7}{2}}\varphi^3\th^2 |v_{s}|^2\mathrm{dx} \mathrm{ds} +C\int_{Q_{b_1,b_2}}(A_2+A_3)\mathrm{dx} \mathrm{ds} \]\\
\ns\ds\q+C\int_{\Sigma_{b_1,b_2}}\(\l^2\mu^3\varphi^2\th^2|\D v|^2+\th^2|\n\D v|^2\)\mathrm{d\Sigma_{b_1,b_2}}  \ea\ee
where $A_j (j=2, 3$) is given by (\ref{A11-12-13});

ii) In the case of $v=\D v=0$ on the boundary, we have
 \bel{0406-a}\ba{ll}\ds
 \mbox{ The left hand side of (\ref{2p14-0}) }\\
 \ns\ds\le C\[||\th\cP v ||^2_{L^2(Q_{b_1, b_2})}+\int_{Q_{b_1,b_2}}\lambda^3\mu^{\frac{7}{2}}\varphi^3\th^2 |v_{s}|^2\mathrm{dx} \mathrm{ds} +C\int_{Q_{b_1,b_2}}(A_2+A_3)\mathrm{dx} \mathrm{ds} \].
 \ea\ee
\el

{\bf Proof of Lemma \ref{lem:5}. } Noting that $\cP v=\b v_{ss}+\a v_s+\D^2 v$, by elementary calculus, we have
\bel{2p13-o}\ba{ll}&\ds
\b \lambda\mu^{\frac{3}{2}}\varphi \th^2\cP v \D v=\b^2 \lambda\mu^{\frac{3}{2}}\varphi\th^2 v_{ss}\D v+\b \lambda\mu^{\frac{3}{2}}\varphi\th^2\D^2 v \D v+\a\b \lambda\mu^{\frac{3}{2}}\varphi\th^2v_s\D v\\
\ns&\ds= \(\b^2 \lambda\mu^{\frac{3}{2}}\varphi\th^2 v_{s}\D v\)_s-\(\b^2 \lambda\mu^{\frac{3}{2}}\varphi\th^2\)_s v_{s}\D v-\n\cdot\(\b^2 \lambda\mu^{\frac{3}{2}}\varphi\th^2 v_{s}\n v_s\)\\
\ns&\ds\q+\frac{\b^2}{2}\n\cdot\(\n \(\lambda\mu^{\frac{3}{2}}\varphi\th^2\)|v_{s}|^2\)-\frac{\b^2}{2}\D \(\lambda\mu^{\frac{3}{2}}\varphi\th^2\) |v_{s}|^2+\b^2 \lambda\mu^{\frac{3}{2}}\varphi\th^2 |\n v_s|^2\\
\ns&\ds\q+\n\cdot\( \b \lambda\mu^{\frac{3}{2}}\varphi\th^2\n\D v \D v\)-\b \lambda\mu^{\frac{3}{2}}\varphi\th^2|\n\D v|^2\\
\ns&\ds\q-\frac{\b}{2}\n\cdot\(\n \(\lambda\mu^{\frac{3}{2}}\varphi\th^2\)|\D v|^2\)+\frac{\b}{2}\D \(\lambda\mu^{\frac{3}{2}}\varphi\th^2\)|\D v|^2+\a\b \lambda\mu^{\frac{3}{2}}\varphi\th^2v_s\D v.
\ea \ee
 Integrating (\ref{2p13-o}) on $Q_{b_1,b_2}$, using the boundary conditions $\ds v=\frac{\pa v}{\pa\nu}=0$ and noting that $\b<0$, we have
\bel{2p14-k}\ba{ll}\ds
\int_{Q_{b_1,b_2}}\(-\b\lambda\mu^{\frac{3}{2}}\varphi \th^2|\n\D v|^2+\b^2 \lambda\mu^{\frac{3}{2}}\varphi \th^2 |\n v_s|^2\)\mathrm{dx} \mathrm{ds} \\
\ns\ds\q+\int_{Q_{b_1,b_2}}\n\cdot\(-\frac{\b}{2} \n \(\lambda\mu^{\frac{3}{2}}\varphi\th^2\)|\D v|^2+\b \lambda\mu^{\frac{3}{2}}\varphi\th^2\n\D v \D v\)\mathrm{dx} \mathrm{ds}\\
\ns\ds\leq C\[||\cP v ||^2_{L^2(Q_{b_1, b_2})}+\int_{Q_{b_1,b_2}}\lambda^3\mu^{\frac{7}{2}}\varphi^3\th^2 |v_{s}|^2\mathrm{dx} \mathrm{ds} +C(A_2+A_3)\].
\ea \ee
Observe that
 \bel{2p14-0l}\ba{ll}\ds
\int_{Q_{b_1,b_2}}\n\cdot\(-\frac{\b}{2} \n \(\lambda\mu^{\frac{3}{2}}\varphi\th^2\)|\D v|^2+\b \lambda\mu^{\frac{3}{2}}\varphi\th^2\n\D v \D v\)\mathrm{dx} \mathrm{ds}\\
\ns\ds \ge-C \int_{\Sigma_{b_1,b_2}}\(\l^2\mu^3\varphi^2|\D v|^2+|\n\D v|^2\)\mathrm{dx} \mathrm{ds}.
\ea\ee
Combining (\ref{2p14-k}) and (\ref{2p14-0l}), we can get (\ref{2p14-0}). Similarly, integrating (\ref{2p13-o}) on $Q_{b_1,b_2}$, using the boundary conditions $\ds v=\D v=0$ and noting that $\b<0$, we can obtain (\ref{0406-a}) immediately. \endpf

 {\bf Proof. }  In \cite[Theorem 1.1]{FLZ}, by taking $\a=0,\  \b=0$ be constants  and $(a^{jk})_{n\t n}=I_n$ (the unit matrix), choosing  $\Psi=-\D\ell$ and $\Phi=0$, a short calculation yields the desired result. \endpf

\section{Proof of Theorem \ref{point-th1}}\label{section2}

In this section, we will give the proof of Theorem \ref{point-th1}. A direct consequence of Lemma \ref{lem:1} and Lemma \ref{lem:02} is the following proposition, which will play an important role in our proof of Theorem \ref{point-th1}.
\bp\label{lem:4}  Under the same conditions of Lemma \ref{lem:1},  we have
  \bel{0112-0a0}\ba{ll}\ds
\cP_1(z)\cP_2(z)+\cP_1(z)\cP_3(z)\\
 \ns\ds
 \ge  \pa_s M+\sum_{j=1}^7\div V_j-\b\Phi z_s^2+\b Z\n^2 \eta\n z_{s}\n z_s-C\sum_{j=1}^4A_j+\sum_{j=1}^4\cH_j\\
 \ns\ds\q+\cB_1|\n z|^2+\cB_2|\D z|^2+\cB_3 z^2+\cB_4|\n\ell\cdot\n z|^2+\mathcal{R}
\ea\ee
where $\cB_j (j=1, 2, 3, 4)$ is given by (\ref{1209-f1}), $\cH_k (k=1, 2, 3, 4)$ is given by (\ref{12101-a}), $M_j$ and  $ V_j$ (j=1,...,7) satisfying respectively (\ref{mv-00a0}), (\ref{m0-a0}) and (\ref{m0v-a0}),   $\mathcal{R}$ is given by \eqref{1210-R4}.
\ep

In view of
the above inequality \eqref{0112-0a0} in Proposition \ref{lem:4},  in order to get the main estimation in Theorem \ref{point-th1}, it is sufficient to estimate the $B_i$ and $H_i$ in \eqref{0112-0a0}.
\bl\label{lem:8} Assume that $\ell\in C^4(Q_{b_1,b_2})$ satisfying Condition \ref{cd1}. Let $\Phi\in C^2(Q_{b_1,b_2};\dbR)$ and $Z\in C^2(Q_{b_1,b_2};\dbR^{n\times n})$ satisfying \eqref{E}.
Then there exists a positive constant $\mu_0$,  such that for any $\mu \ge \mu_0$, one can find two positive constants $\l_0=\l_0(\mu)$, $c_1>0$ and $C = C(\mu)$,  such that \bel{1215-aa0}\ba{ll}\ds
-\b\Phi z_s^2+\cB_1|\n z|^2+\cB_2|\D z|^2+\cB_3 z^2+\cB_4|\n\ell\cdot\n z|^2+\sum_{j=1}^4\cH_j\\
 \ns\ds
 \ge \n\cdot V_9+\b^2\lambda^3\mu^{\frac{7}{2}}\varphi^3|\n\eta|^4 z_s^2+2\l^3\mu^4\varphi^3|\n\eta|^4\th^2 |\D v |^2+3\lambda^3\mu^4\varphi^3|\n\eta|^4|\n\ell\cdot\n z|^2\\
\ns\ds\q+3\lambda\mu^2\varphi|\n\eta\cdot\n\D z|^2+\b Z\n^2\eta\n z_s\n z_s+64\lambda\mu\varphi\n^2\eta(\n^2 z\n\ell)(\n^2 z\n\ell)-C\sum_{j=1}^4A_j,
\ea\ee
where  $V_9$ is given by \eqref{mv-0331a}.
\el
{\bf Proof of Lemma \ref{lem:8} } The proof relies on careful study of
the lower bounds of each terms of L.H.S. in \eqref{0112-0a0}.
 We divide the proof into several steps.

\ms

{\bf Step 1:  Estimations of $\cH_j$ ($j=1, 2, 3, 4)$. }

Noting that $\ds  Z=8\l\mu\varphi$ and recalling \eqref{12101-a} the definitions of $\cH_i$ and condition \eqref{cd1}, we have
\bel{12101-a1}\ba{ll}\ds
  \cH_1=8\l\mu^2\varphi|\n\eta\cdot\n\D z|^2,\q \cH_2 =0.
 \ea\ee

Further, by (\ref{12101-a}) and \eqref{0421-a}, we know that
 \bel{12101-a3en1}\ba{ll}\ds
  \cH_3&\ds
  =32\lambda^3\mu^4\varphi^3|\n^2 z\n\eta\n\eta|^2+32\lambda\mu\varphi\n^2\eta(\n^2 z\n\ell)(\n^2 z\n\ell)\\
  \ns&\ds\q+32\lambda\mu\varphi\(\n^2 z:\n^2\eta\)\n^2 z\n\ell\n\ell\\
  \ea\ee
  By using divergence formula and recalling \eqref{A11-12-13}, we get
  \bel{12101-a3en}\ba{ll}\ds
  \cH_3&\ge 32\n\cdot\[\lambda\mu\varphi\(\n^2 z:\n^2\eta\)\(\n z\cdot\n\ell\)\n\ell-\lambda\mu\varphi\(\n z\cdot\n\ell\)\n^2\eta\(\n^2 z\n\ell\)\]\\
 \ns&\ds\q+64\lambda\mu\varphi\n^2\eta\(\n^2 z\n\ell\)\(\n^2 z\n\ell\)-CA_1-CA_2
  \ea\ee
 Similarly, we have
 \bel{12101-a4}\ba{ll}\ds
 \cH_4=-16A\lambda\mu^2\varphi\n^2 z\n\eta\n\eta\D z\\
 \ns\ds\qq\ge -16\n\cdot\(A\lambda\mu^2\varphi\(\n z\cdot\n\eta\)\D z\n\eta\)+16\lambda^3\mu^4\varphi^3|\n\eta |^2 \(\n\D z\cdot\n\eta\)\cdot\(\n z\cdot\n\eta\)\\
 \ns\ds\qq\q-CA_1-CA_2+ 16\lambda^3\mu^5\varphi^3  |\n\eta|^4 \D z\(\n z\cdot\n\eta\)\\
 \ns\ds\qq\ge -4\lambda\mu^2\varphi |\n\D z\cdot\n\eta|^2-20\lambda^5\mu^6\varphi^5 |\n\eta |^4 |\n z\cdot\n\eta |^2\\
\ns\ds\qq\q+\lambda^{-1}\varphi^{-1}\(2\lambda\mu\varphi \n\D z\cdot\n\eta+4\lambda^3\mu^3\varphi^3 |\n\eta |^2 \n z\cdot\n\eta\)^2\\
\ns\ds\qq\q-16\n\cdot\(A\lambda\mu^2\varphi\(\n z\cdot\n\eta\)\D z\n\eta\)-CA_1-CA_2-16\lambda\mu^4\varphi|\n\eta |^4 |\D z|^2.
\ea\ee

By (\ref{12101-a1})--(\ref{12101-a4}), we have
 \bel{0110-0b0}\ba{ll}&\ds
 \cH_1+\cH_2+\cH_3+\cH_4\\
 \ns&\ds\ge 32\n\cdot\[\lambda\mu\varphi\n^2 z:\n^2\eta\n z\cdot\n\ell\n\ell-\lambda\mu\varphi\n^2\eta\(\n z\cdot\n\ell\)\(\n^2 z\n\ell\)\]\\
 \ns&\ds\q-16\n\cdot\(A\lambda\mu^2\varphi\n z\cdot\n\eta\D z\n\eta\)-20\lambda^5\mu^6\varphi^5 |\n\eta |^4 |\n z\cdot\n\eta |^2+4\l\mu^2\varphi|\n\eta\cdot\n\D z|^2\\
 \ns&\ds\q+64\lambda\mu\varphi\n^2\eta(\n^2 z\n\ell)(\n^2 z\n\ell)-CA_1-CA_2-16\lambda\mu^4\varphi|\n\eta |^4 |\D z|^2.
 \ea\ee

\ms

{\bf Step 2: Estimations of $\cB_j$ ($j=1,\cdots, 4)$. }

Noting that  $A=|\n \ell|^2$ and $\ds
 \Phi=-\b\lambda^3\mu^{\frac{7}{2}}\varphi^3|\n\eta|^4 $, a short calculation shows that
 \bel{1211-b}\left\{\ba{ll}\ds
 \cB_1=-8A\n A\cdot\n\ell-2A\D^2\ell+2A\b\lambda^3\mu^{\frac{7}{2}}\varphi^3|\n\eta|^4,\\
 \ns\ds\cB_2=4\n A\cdot\n\ell+\D^2\ell-\b\lambda^3\mu^{\frac{7}{2}}\varphi^3|\n\eta|^4,\\
 \ns\ds\cB_3=4A^2\n A\cdot\n\ell+A^2\D^2\ell-A^2\b\lambda^3\mu^{\frac{7}{2}}\varphi^3|\n\eta|^4,\\
 \ns\ds\cB_4=8A\lambda\mu^2\varphi|\n\eta|^2+16\n\cdot\(A\n\ell\)-16A\D\ell-4\D^2\ell+4\b\lambda^3\mu^{\frac{7}{2}}\varphi^3|\n\eta|^4.
 \ea\right.\ee

  On the one hand, by (\ref{1211-b}), we know that
 \bel{r11a1}\ba{ll}\ds
\cB_1|\n z|^2+\cB_2|\D z|^2+\cB_3 z^2\\
\ns\ds=A^2\(4\n A\cdot\n\ell+\D^2\ell-\b\lambda^3\mu^{\frac{7}{2}}\varphi^3|\n\eta|^4\) | z |^2-2A\(4\n A\cdot\n\ell+\D^2\ell-\b\lambda^3\mu^{\frac{7}{2}}\varphi^3|\n\eta|^4\)  |\n z|^2\\
\ns\ds\q+\(4\n A\cdot\n\ell+\D^2\ell-\b\lambda^3\mu^{\frac{7}{2}}\varphi^3|\n\eta|^4\)|\D z |^2\\
\ns\ds=\(4\n A\cdot\n\ell+\D^2\ell-\b\lambda^3\mu^{\frac{7}{2}}\varphi^3|\n\eta|^4\)\(A z+\D z \)^2\\
\ns\ds\q-2A\(4\n A\cdot\n\ell+\D^2\ell-\b\lambda^3\mu^{\frac{7}{2}}\varphi^3|\n\eta|^4\)   z\D z\\
\ns\ds\q- 2A\(4\n A\cdot\n\ell+\D^2\ell-\b\lambda^3\mu^{\frac{7}{2}}\varphi^3|\n\eta|^4\)   |\n z|^2\\
\ea\ee
Using divergence formula and recalling \eqref{A11-12-13}, we have
\bel{r11a}\ba{ll}\ds
-2A\(4\n A\cdot\n\ell+\D^2\ell-\b\lambda^3\mu^{\frac{7}{2}}\varphi^3|\n\eta|^4\)   z\D z
 - 2A\(4\n A\cdot\n\ell+\D^2\ell-\b\lambda^3\mu^{\frac{7}{2}}\varphi^3|\n\eta|^4\)   |\n z|^2\\
\ns\ds\ge -8\n\cdot\(A\n A\cdot\n\ell z\n z\)
+2\n\cdot\(A\b\lambda^3\mu^{\frac{7}{2}}\varphi^3|\n\eta|^4 z\n z\)+8\n\(A\n A\cdot\n\ell\)\cdot\n z   z\\
\ns\ds\qq-2\n\(A\b\lambda^3\mu^{\frac{7}{2}}\varphi^3|\n\eta|^4\)\cdot\n z z- CA_1-CA_2 \\
\ea\ee
By using \eqref{0421-a} and \eqref{A11-12-13}, we get that
\bel{r11a2-0}\ba{ll}\ds
8\n\(A\n A\cdot\n\ell\)\cdot\n z   z=\(\n A\cdot \n \ell\)\(\n A\cdot\n z\) z+A\n^2A\n \ell\n z z+A\n^2\ell\n A\n z z\\
\ns\ds \geq -8\lambda^5\mu^7\varphi^5|\n \eta|^6|\n \eta\cdot \n z|z- CA_1\geq -4\lambda^4\mu^7\varphi^4|\n \eta|^4|\n z\cdot\n \eta|^2   - CA_1.\\
\ea\ee
Similarly,
\bel{r11a2}\ba{ll}\ds
-2\n\(A\b\lambda^3\mu^{\frac{7}{2}}\varphi^3|\n\eta|^4\)\cdot\n z z \geq - CA_1.\\
\ea\ee

On the other hand, by (\ref{0421-b}) and inequality $(a+b+c)^2\leq 3(a^2+b^2+c^2)$, we have
\bel{r13a}\ba{ll}\ds
(\n A \cdot \n\ell) \th^2 |\D v |^2=(\n A \cdot \n\ell)\(\D z+A z-2\n\ell\cdot\n z-\D\ell z\)^2\\
\ns\ds\le 3\(\n A\cdot\n\ell\)\(A z+\D z \)^2 +12\(\n A\cdot\n\ell\) | \n\ell\cdot\n z |^2+CA_1 \\
\ns\ds\le 3\(\n A\cdot\n\ell\)\(A z+\D z \)^2 +12\lambda^5\mu^{6}\varphi^5|\n\eta|^4 | \n\eta\cdot\n z |^2+CA_1.
\ea\ee
  By \eqref{A11-12-13} and \eqref{1211-b}, we have
\bel{b-o4}
\cB_4|\n\ell\cdot \n z|^2\ge 36\lambda^5\mu^6\varphi^5|\n\eta|^4| \n\eta\cdot\n z |^2-A_1,
 \ee
 this, together with \eqref{r11a1}-\eqref{b-o4}, we have
\bel{b-4}\ba{ll}\ds
\cB_1|\n z|^2+\cB_2|\D z|^2+\cB_3 z^2+\cB_4|\n\ell\cdot\n z|^2+C(A_1+A_2) \\
\ns\ds \geq (24\lambda\varphi-4\mu)\lambda^4\mu^6\varphi^4|\n\eta|^4| \n\eta\cdot\n z |^2+\(\n A\cdot\n\ell+\D^2\ell-\b\lambda^3\mu^{\frac{7}{2}}\varphi^3|\n\eta|^4\)\(A z+\D z \)^2\\
\ns\ds\q + \n\cdot\(-8A\n A\cdot\n\ell z\n z
+2A\b\lambda^3\mu^{\frac{7}{2}}\varphi^3|\n\eta|^4 z\n z\)
\ea\ee

{\bf Step 3:  Estimations of remainder term $\mathcal{R}$ } Taking \eqref{E} into \eqref{1210-R4}, by \eqref{0421-a}, \eqref{0227-01} and \eqref{A11-12-13}, we have
\bel{}\ba{ll}\ds
|\mathcal{R}|\leq C\[\lambda^{3}\mu^{\frac{9}{2}}\varphi^3|\n \eta\cdot\n\D z||z|+\lambda^{3}\mu^{\frac{7}{2}}\varphi^3|\n\D z||z|+\lambda^{3}\mu^{\frac{9}{2}}\varphi^3|\n\eta\cdot\n z||\D z| \\
\ns\ds\qq\q +\lambda^{3}\mu^{\frac{7}{2}}\varphi^3|\n z||\D z|
+\lambda^{5}\mu^{\frac{13}{2}}\varphi^5|\n\eta\cdot\n z|| z|+\lambda^{5}\mu^{\frac{11}{2}}\varphi^5|\n z|| z|+\lambda^{3}\mu^{4}\varphi^3|\n\eta\cdot\n z|| \n^2 z|
 \\
\ns\ds\qq\q+\lambda^{3}\mu^{3}\varphi^3|\n z||\n^2 z|
+\lambda\mu^{2}\varphi|\n\eta\cdot\n \D z||\n^2 z|+ \lambda\mu^{}\varphi|\n \D z||\n^2 z|
\].
\ea\ee
In view of \eqref{A11-12-13}, we have
\bel{R1}
\mathcal{R}\geq -C (A_1+A_2+A_4).
\ee

{\bf Step 4:  End of the proof. }
By (\ref{0421-a}) and (\ref{0421-02}), we conclude that there is a $\mu_0>0$, for any $\mu\ge\mu_0$ and $\l>C\mu$, we have
 \bel{0108-0a}\left\{\ba{ll}\ds
16\lambda\mu^4\varphi|\n\eta |^4 |\D z|^2\leq  CA_1+CA_2, \\
\ns\ds (24\lambda\varphi-4\mu)\lambda^4\mu^6\varphi^4\geq 23 \lambda^5\mu^6\varphi^5\\
\ns\ds \n A\cdot\n\ell+\D^2\ell-\b\lambda^3\mu^{\frac{7}{2}}\varphi^3|\n\eta|^4\ge 0,\\
\ns\ds -\b\Phi \geq \frac{1}{2}\beta^2\lambda^3\mu^{\frac{7}{2}}\varphi^3|\n\eta|^4.
\ea\right.\ee
  In view of \eqref{0112-0a0} in Proposition \ref{lem:4},
combining \eqref{0110-0b0}, (\ref{b-4}), \eqref{R1} and (\ref{0108-0a}), we obtain \eqref{1215-aa0}. This concludes the proof. \endpf

\ms

Now we are in a position to prove Theorem \ref{point-th1}.
It is only need to  get back the lower order estimate of $v$ by using $\D v$.

 {\bf Proof of Theorem \ref{point-th1}. }
 Set
  $$
  y=\varphi^{3\over 2}|\n\eta|^2 v,\qq u=\th y.
  $$
Then
  \bel{v-y}
  \th\D y=e^{\lambda\varphi}\varphi^{3\over 2}|\n\eta|^2 \D v+2e^{\lambda\varphi}\n (\varphi^{3\over 2}|\n\eta|^2)\cdot\n v+e^{\lambda\varphi}\D(\varphi^{3\over 2}|\n\eta|^2) v.
 \ee
Applying Lemma \ref{2l1} with \eqref{v-y}, we have
   \bel{2aao2}\ba{ll}
\th^2\varphi^3|\n\eta|^4| \D v|^2+\n\cdot V_8&\ds\ge
 c_0\[\lambda^3\mu^4\varphi^6\th^2|\nabla\eta|^8|v|^2 +\lambda\mu^2\varphi^4\th^2|\nabla\eta|^6|\n v|^2\]\\
 \ns&\ds\q- C\th^2\[(\lambda^3\mu^3\varphi^6+\l^2\mu^4\varphi^5)|v|^2+\l\mu\varphi|\n v|^2\],
	 \ea\ee
where  $V_8$ is given by \eqref{mv-0331a}.
Multiplying both sides of (\ref{2aao2}) by $\lambda^3\mu^4$, recalling \eqref{A11-12-13} again, we obtain
   \bel{v-1}\ba{ll}\ds
 \lambda^3\mu^4\varphi^3 \th^2|\n\eta |^4 |\D v |^2+ \lambda^3\mu^4\n\cdot V_8\\
 \ns\ds\ge  c_0\[\lambda^6\mu^8\varphi^6\th^2|\nabla\eta|^8|v|^2+\lambda^4\mu^6\varphi^4\th^2|\nabla\eta|^6|\n v|^2\]- C(A_1+A_2).
  \ea\ee

Finally, taking  (\ref{v-1}) into \eqref{1215-aa0} in Lemma \ref{lem:8}, recalling \eqref{Aii},  we can get the desired result immediately.
\endpf

 \ms

\section{ Proof of Theorem \ref{thm:2}}\label{section4-0}

In this section, we will use Theorem \ref{point-th1}  to prove  Theorem \ref{thm:2}.

{\bf Proof of Theorem \ref{thm:2}. }
 In view of (\ref{point-th1-01}) in Theorem \ref{point-th1},  integrating on $Q_{b_1, b_2}$,  we have
 \bel{0112-ax0}\ba{ll}\ds
 c\int_{Q_{b_1, b_2}}\th^2\(\lambda^6\mu^8\varphi^6|\n\eta |^8 v^2+\lambda^4\mu^6\varphi^4 |\n\eta |^6 |\n v |^2+\l^3\mu^4\varphi^3|\n\eta|^4|\D v |^2\)\mathrm{dxds}\\
  \ns\ds\q +c\int_{Q_{b_1, b_2}}\(\l^3\mu^{\frac{7}{2}}\varphi^3|\n\eta|^4z_s^2+ \lambda\mu^2\varphi|\n\eta\cdot\n\D z|^2+\lambda^5\mu^6\varphi^5|\n\eta|^4|\n\eta\cdot\n z|^2\)\mathrm{dxds}\\
  \ns\ds\q+8\b\int_{Q_{b_1, b_2}}\lambda\mu\varphi\n^2\eta \n z_s\n z_s\mathrm{dx}\mathrm{ds}+64\int_{Q_{b_1, b_2}}\lambda\mu\varphi\n^2\eta(\n^2 z\n\ell)(\n^2 z\n\ell)\mathrm{dxds}\\
  \ns\ds\q +\int_{Q_{b_1, b_2}}\(\pa_s M+ \n \cdot V\)\mathrm{dxds}\\
 \ns\ds\le C\int_{Q_{b_1, b_2}}\th^2\( |\cP v|^2+\l^4\mu^5\varphi^4|\n v|^2+\l^6\mu^7\varphi^6v^2+\l^2\mu^3\varphi^2|\n^2 v|^2\\
 \ns\ds\qq\qq\qq\qq\qq\qq\qq\q+\b^2\l^2\mu^3\varphi^2v_s^2+\mu|\n\D v|^2\)\mathrm{dxds}\].
 \ea\ee
In order to get \eqref{i1i2-c1}, we divide our proof into four steps.
\ss

{\bf Step 1. Estimation of boundary terms. } Noting that $z=\th v$,  by (\ref{0305-a}) and (\ref{mv-00a0}),  we have
 \bel{0112-0v}
 \int_{Q_{b_1, b_2}} \pa_s M \mathrm{dxds}=0.
 \ee
  Since $\eta$ satisfies (\ref{1224-0a}), by (\ref{m0-a0})-(\ref{mv-0331a}),  we have
 \bel{0112-k0}
 \int_{Q_{b_1, b_2}} \n\cdot V_8\mathrm{dxds}=0,\q \int_{Q_{b_1, b_2}} \n\cdot (V_1+V_2+V_3+V_9) \mathrm{dxds}=0.
 \ee
 By Stoke's formula,
 \bel{0112-0d0}\ba{ll}\ds
 \int_{Q_{b_1, b_2}} \n\cdot(V_4+V_5+V_6+V_7) \mathrm{dxds}= \int_{\Sigma_{b_1,b_2}} \cH\cdot\nu \mathrm{d\Sigma_{b_1,b_2}}
 \ea\ee
   where $\mathrm{d\Sigma_{b_1,b_2}}=\mathrm{d\pa\O ds}$ and
   \bel{mv-0}\ba{ll}\ds
  \cH=-4\n\D z(\n\ell\cdot\n\D z)+2\n\ell|\n\D z|^2-4\n\D z(\n^2\ell:\n^2 z)\\
  \ns\ds\qq+4\n^2 \ell\cdot(\n^2 z\cdot\n\D z)-4\n^2 z\cdot(\n^2\ell\cdot\n\D z)-2\D\ell\D z\n\D z\\
  \ns\ds\qq+\n\D\ell|\D z|^2+4A \n^2 z\n\ell\n^2 z-2A\n\ell |\n^2 z |^2-4A\n\ell|\D z|^2\\
  \ns\ds\qq-16\(\n^2 z\n\ell\n\ell\)\n^2 z\n\ell+8| \n^2 z\n\ell |^2\n\ell+( Z\n^2 z:\n^2\eta)\n\D z\\
\ns\ds\qq-Z\n^2\eta\( \n^2  z\cdot\n\D z\) +Z\n^2 z\( \n^2 \eta\cdot\n\D z\)-4\n^2 z \n\D\ell\D z.
  \ea\ee

Next, we estimate $\cH\cdot\nu$. We can decompose $\n$ into  the normal component $\frac{\pa}{\pa\nu}$ and the tangential component $\frac{\pa}{\pa\tau}$, such that by \eqref{1224-0a},
 \bel{1225-0b}
\q \n\eta=\frac{\pa \eta}{\pa\nu}\nu+\frac{\pa \eta}{\pa\tau}\tau=\frac{\pa \eta}{\pa\nu}\nu \q \mbox{ and } \q\frac{\pa\eta}{\pa\nu}< 0, \q \mbox{on}~ \pa\O.
 \ee
 Moreover, by $z=\th v$ and the boundary conditions satisfied by $v$, we have
 \bel{bound:z}
z=0, z_s=0, \n z=0, \frac{\pa^2 z}{\pa\tau\pa \nu}=0.
 \ee
  For simplicity, we define
 \bel{b24-ae}\ba{ll}\ds
A_5= \lambda^2\mu^2\int_{\Sigma_{b_1,b_2}}\varphi^2\(\lambda\varphi+\mu\)|\D z|^2\mathrm{d\Sigma_{b_1,b_2}}+\int_{\Sigma_{b_1,b_2}}\(\lambda\varphi+\mu\)|\frac{\pa\D z}{\pa\nu}|^2\mathrm{d\Sigma_{b_1,b_2}}\\
\ns\ds\qq+\lambda^2\mu^3\int_{\Sigma_{b_1,b_2}}\varphi^2  |\n^2 z|^2\mathrm{d\Sigma_{b_1,b_2}}+\lambda^{\frac{1}{2}}\mu\int_{\Sigma_{b_1,b_2}} \varphi^{\frac{1}{2}} |\n\D z|^2\mathrm{d\Sigma_{b_1,b_2}}.
\ea\ee
Then recalling \eqref{0421-a}, we get
 \bel{1225-c0}\ba{ll}\ds
 \int_{\Sigma_{b_1,b_2}}\cH\cdot\nu\mathrm{d\Sigma_{b_1,b_2}}\\
 \ns\ds\ge 12\lambda^3\mu^3\int_{\Sigma_{b_1,b_2}}\varphi^3\(\frac{\pa\eta}{\pa\nu}\)^3 |\n^2 z\cdot\nu|^2\mathrm{d\Sigma_{b_1,b_2}} +2\lambda\mu\int_{\Sigma_{b_1,b_2}}\varphi\frac{\pa\eta}{\pa\nu}|\n\D z|^2\mathrm{d\Sigma_{b_1,b_2}}\\
 \ns\ds\q-4\lambda\mu\int_{\Sigma_{b_1,b_2}}\varphi\frac{\pa\eta}{\pa\nu}|\frac{\pa\D z}{\pa\nu}|^2\mathrm{d\Sigma_{b_1,b_2}}-2\lambda^3\mu^3\int_{\Sigma_{b_1,b_2}}\varphi^3\(\frac{\pa\eta}{\pa\nu}\)^3\|\n^2 z|^2 \mathrm{d\Sigma_{b_1,b_2}} \\
 \ns\ds\q-4\lambda^3\mu^3\int_{\Sigma_{b_1,b_2}}\varphi^3\(\frac{\pa\eta}{\pa\nu}\)^3 |\D z|^2\mathrm{d\Sigma_{b_1,b_2}}-16\lambda^3\mu^3\int_{\Sigma_{b_1,b_2}}\varphi^3\(\frac{\pa\eta}{\pa\nu}\)^3 |\frac{\pa^2 z}{\pa\nu^2}|^2\mathrm{d\Sigma_{b_1,b_2}}-CA_5.
\ea\ee

\ss
Since $\frac{\pa\eta}{\pa\nu}< 0$ on $\pa \Omega$, we now only need to deal with the first two terms in the R.H.S of  \eqref{1225-c0}.
   By \eqref{bound:z}, we have
\bel{b24-cei}\ba{ll}\ds
12\int_{\Sigma_{b_1,b_2}}\lambda^3\mu^3\varphi^3\(\frac{\pa\eta}{\pa\nu}\)^3 |\n^2 z\cdot\nu|^2\mathrm{d\Sigma_{b_1,b_2}}&\ds=12\int_{\Sigma_{b_1,b_2}}\lambda^3\mu^3\varphi^3\(\frac{\pa\eta}{\pa\nu}\)^3 |\frac{\pa^2 z}{\pa\nu^2}+\frac{\pa^2 z}{\pa\nu\pa\tau}|^2\mathrm{d\Sigma_{b_1,b_2}}\\
\ns&\ds=12\int_{\Sigma_{b_1,b_2}}\lambda^3\mu^3\varphi^3\(\frac{\pa\eta}{\pa\nu}\)^3 |\frac{\pa^2 z}{\pa\nu^2}|^2\mathrm{d\Sigma_{b_1,b_2}},
\ea\ee
and
\bel{eq:4.10}\ba{ll}\ds
2\int_{\Sigma_{b_1,b_2}}\lambda\mu\varphi\frac{\pa\eta}{\pa\nu}|\n\D z|^2\mathrm{d\Sigma_{b_1,b_2}}\\
\ns\ds=2\int_{\Sigma_{b_1,b_2}}\lambda\mu\varphi\frac{\pa\eta}{\pa\nu}|\frac{\pa\D z}{\pa\nu}|^2\mathrm{d\Sigma_{b_1,b_2}}+2\int_{\Sigma_{b_1,b_2}}\lambda\mu\varphi\frac{\pa\eta}{\pa\nu}|\frac{\pa\D z}{\pa\tau}|^2\mathrm{d\Sigma_{b_1,b_2}}.
\ea\ee
The first one is  obviously bounded by $\int_{\Sigma_{b_1,b_2}}-16\lambda^3\mu^3\varphi^3\(\frac{\pa\eta}{\pa\nu}\)^3 |\frac{\pa^2 z}{\pa\nu^2}|^2\mathrm{d\Sigma_{b_1,b_2}}$. So our task now is to estimate $2\int_{b_1}^{b_2}\int_{\pa\O}\lambda\mu\varphi\frac{\pa\eta}{\pa\nu}|\frac{\pa\D z}{\pa\tau}|^2\mathrm{d\Sigma_{b_1,b_2}}$.
In view of \eqref{0421-a}, by \eqref{bound:z}, we have
 \bel{eq:4.11}
 2\int_{\Sigma_{b_1,b_2}}\lambda\mu\varphi\frac{\pa\eta}{\pa\nu}|\frac{\pa\D z}{\pa\tau}|^2\mathrm{d\Sigma_{b_1,b_2}}=2\int_{\Sigma_{b_1,b_2}}\lambda\mu\varphi\th^2\frac{\pa\eta}{\pa\nu}|\frac{\pa\D v}{\pa\tau}|^2\mathrm{d\Sigma_{b_1,b_2}}.
 \ee
 Notice that for any $x\in\pa \O$, \bel{th0}\th(s,x)=\th_0(s)=e^{\l\xi(s,x_1)},\q x_1\in \pa\O,\ee since $\xi$ is independent with space variable $x$ on the boundary $\pa\O$.
 Moreover, there exist $c,C>0$, such that
 \bel{th1}
  c\th \leq \th_0\leq C\th,\q  |(\th_0(s))_s|\leq C\l\mu\varphi\th, \q |(\th_0(s))_{ss}|\leq C\l^2\mu^2\varphi^2\th.
 \ee

 Then we have
 \bel{b24-de}\ba{ll}\ds
2\int_{\Sigma_{b_1,b_2}}\lambda\mu\varphi\th^2\frac{\pa\eta}{\pa\nu}|\frac{\pa\D v}{\pa\tau}|^2\mathrm{d\Sigma_{b_1,b_2}}
\ge-C\int_{b_1}^{b_2}\lambda\mu\varphi \th_0^2(s)||\D v||^2_{H^1(\pa\Omega)}\mathrm{ds}.
\ea\ee
 Recalling the fact that $\frac{\pa\eta}{\pa\nu}<0$ on $\pa\O$ again, by using the interpolation inequality (see \cite{AdamsBook}) and  Young's inequality, we deduce
\bel{b24-de}\ba{ll}\ds
2\int_{\Sigma_{b_1,b_2}}\lambda\mu\varphi\th^2\frac{\pa\eta}{\pa\nu}|\frac{\pa\D v}{\pa\tau}|^2\mathrm{d\Sigma_{b_1,b_2}}\\
\ns\ds\ge -C\int_{b_1}^{b_2}\lambda\mu\varphi \th_0^2(s)||\D v||^{\frac{2}{3}}_{L^2(\pa\Omega)} ||\D v||^{\frac{4}{3}}_{H^{\frac{3}{2}}(\pa\Omega)}\mathrm{d\Sigma_{b_1,b_2}}\\
\ns\ds\ge \int_{\Sigma_{b_1,b_2}}\lambda^3\mu^3\varphi^3\(\frac{\pa\eta}{\pa\nu}\)^3 \th_0^2(s)|\D v|^2\mathrm{d\Sigma_{b_1,b_2}} -C\int_{b_1}^{b_2} \th_0^2(s)||\D v||^2_{H^{\frac{3}{2}}(\pa\Omega)}\mathrm{ds}.
\ea\ee
Observing that $\th_0(s)\D v=\D z$ on $\pa \O$ and  using  the continuity of the trace operator, we have
\bel{eq:4.14}\ba{ll}\ds
2\int_{\Sigma_{b_1,b_2}}\lambda\mu\varphi\th^2\frac{\pa\eta}{\pa\nu}|\frac{\pa\D v}{\pa\tau}|^2\mathrm{d\Sigma_{b_1,b_2}}\\
\ns\ds\ge\int_{\Sigma_{b_1,b_2}} \lambda^3\mu^3\varphi^3\(\frac{\pa\eta}{\pa\nu}\)^3|\D z|^2\mathrm{d\Sigma_{b_1,b_2}} -C\int_{b_1}^{b_2}\th^2_0(s)||v||^2_{H^4(\Omega)}\mathrm{ds}.
\ea\ee

 According to  \cite[Theorem 2.20]{ggs}(or \cite[Theorem 15.2]{adn}), we can claim that $\ds  ||v||^2_{H^4(\Omega)}$ can be bounded by $||\D^2 v||_{L^2(\O)}$ and some lower order terms. More precisely,  by taking $\ds m=2,\  k=4, \ p=2, \ a_\b=0, \b_{1, \a}=0,\ \b=0,\ \f_\b=f,\ \b_{2,\a}=1$ and $m_1=0, m_2=0, h_0=h_1=0,  u=v$ in  \cite[Theorem 2.20]{ggs}, then  the solution $v$ satisfies
  \bel{0324}\ba{ll}\ds
  ||v||_{H^4(\O)}\le C||\D^2 v||_{L^2(\O)}.
  \ea\ee
  Therefore, combing \eqref{1225-c0}, \eqref{b24-cei},\eqref{eq:4.10},\eqref{eq:4.11},\eqref{eq:4.14} and \eqref{0324}, we have
 \bel{b24-a3}\ba{ll}\ds
 \int_{\Sigma_{b_1,b_2}} \cH\cdot \nu \mathrm{d\Sigma_{b_1,b_2}}\\
\ns\ds\ge-2\int_{\Sigma_{b_1,b_2}}\lambda\mu\varphi\frac{\pa\eta}{\pa\nu}|\frac{\pa\D z}{\pa\nu}|^2 \mathrm{d\Sigma_{b_1,b_2}} -2\int_{b_1}^{b_2}\int_{\pa\O}\lambda^3\mu^3\varphi^3\(\frac{\pa\eta}{\pa\nu}\)^3\|\n^2 z|^2\mathrm{d\Sigma_{b_1,b_2}}\\
 \ns\ds\q-3\int_{\Sigma_{b_1,b_2}}\lambda^3\mu^3\varphi^3\(\frac{\pa\eta}{\pa\nu}\)^3 |\D z|^2\mathrm{d\Sigma_{b_1,b_2}}-4\int_{b_1}^{b_2}\int_{\pa\O}\lambda^3\mu^3\varphi^3\(\frac{\pa\eta}{\pa\nu}\)^3 |\frac{\pa^2 z}{\pa\nu^2}|^2\mathrm{d\Sigma_{b_1,b_2}}\\
\ns\ds\q-CA_5-C\int_{b_1}^{b_2}\th_0^2(s)||\D^2 v||^2_{L^2(\Omega)}\mathrm{ds}.
\ea\ee
Then there are  $\mu_0>0$ and  $\l_0>0$, for all $\l\ge\l_0$, there is  a $c_1>0$ such that
  \bel{1226-h0}\ba{ll}\ds
  \int_{Q_{b_1, b_2}}\cH\cdot\nu \mathrm{dxds}\\
  \ns\ds\ge c_1\int_{\Sigma_{b_1,b_2}}\th^2\(\lambda\mu\varphi|\n\D v|^2+\l^3\mu^3\varphi^3|\D v|^2 \)\mathrm{d\Sigma_{b_1,b_2}}-C\int_{b_1}^{b_2}\th_0^2(s)||\D^2 v||^2_{L^2(\Omega)}\mathrm{ds}.
  \ea\ee

By elementary calculus, we know that
  \bel{1225-ve0}\ba{ll}\ds
  \th_0^2(s)\(|\b v_{ss}|^2+|\D^2 v|\)\\
  \ns\ds=\th_0^2(s)|\cP v-\a v_s|^2-2\b \n\cdot \(\th_0^2(s)(v_{ss}\n\D v-\n v_{ss}\D v)\)\\
  \ns\ds\q-2\b\(\th_0^2(s)\D v_s \D v\)_s+2\b \th_0^2(s)|\D v_s|^2\\
  \ns\ds\q+\b \(\th_0^2(s)\)_s|\D v|^2\)_s-\b\(\th_0^2(s)\)_{ss}|\D v|^2.
    \ea\ee
Noting that $v=0$ on $\partial \Omega $,  $ v(b_1, x)= v(b_2, x)=0$ and $\b<0$, integrating (\ref{1225-ve0}) on $(b_1, b_2)\t\O$, we have
 \bel{1225-vf0}\ba{ll}\ds
   \int_{Q_{b_1, b_2}} \th_0^2(s)|\D^2 v|\mathrm{dxds}\le  \int_{Q_{b_1, b_2}}  \th_0^2(s)|\cP v-\a v_s|^2\mathrm{dxds}+C  \int_{Q_{b_1, b_2}} |\b|\|\(\th_0^2(s)\)_{ss}\||\D v|^2\mathrm{dxds}.
    \ea\ee
Combining  (\ref{th1}), (\ref{1226-h0}) and (\ref{1225-vf0}), we conclude that
 \bel{1225}\ba{ll}\ds
   \int_{\Sigma_{b_1,b_2}} \cH\cdot \nu \mathrm{d\Sigma_{b_1,b_2}}
    \ge c_1\int_{\Sigma_{b_1,b_2}}\th^2\(\lambda\mu\varphi|\n\D v|^2+\l^3\mu^3\varphi^3|\D v|^2 \)\mathrm{d\Sigma_{b_1,b_2}} \\ \ns\ds \qq\qq\qq\qq \qq  -C\int_{Q_{b_1, b_2}}  \th^2|\cP v-\a v_s|^2\mathrm{dxds}-C  \int_{Q_{b_1, b_2}} \l^2\mu^2\varphi^2\th^2|\D v|^2\mathrm{dxds}.
    \ea\ee

  \ms

{\bf Step 2. Estimation of $64\int_{Q_{b_1, b_2}}\lambda\mu\varphi\n^2\eta(\n^2 z\n\ell)(\n^2 z\n\ell)\mathrm{dxds}$.} By using the boundary condition $\ds\frac{\pa z}{\pa\nu}=0$, we have
 \bel{i1i548-03}\ba{ll}\ds
 64 \int_{Q_{b_1, b_2}}\lambda\mu\varphi\n^2\eta\(\n^2 z\n\ell\)\(\n^2 z\n\ell\)\mathrm{dxds}\\
  \ns\ds\ge-C \int_{Q_{b_1, b_2}}\lambda^3\mu^3\varphi^3|\n^2 z\n\eta|^2\mathrm{dxds}\\
  \ns\ds\ge C \int_{Q_{b_1, b_2}}\lambda^3\mu^3\varphi^3\(\n\D z\cdot\n\eta\)\(\n z\cdot\n\eta\)\mathrm{dx}\mathrm{ds}-C\int_{Q_{b_1, b_2}}\(A_1+A_2\)\mathrm{dxds}\\
 \ns\ds\ge -C \int_{Q_{b_1, b_2}}\(\lambda\mu\varphi|\n\D z\cdot\n\eta|^2+\lambda^5\mu^5\varphi^5|\n z\cdot\n\eta|^2+A_1+A_2\)\mathrm{dxds}.
 \ea\ee

{\bf Step 3. Estimation of $8\b\int_{Q_{b_1, b_2}}\lambda\mu\varphi\n^2\eta \n z_s\n z_s\mathrm{dx}\mathrm{ds}$.} Recalling $z=\th v$ and condition \eqref{cd1} satisfied by $\th$,
\bel{}\ba{ll}\ds
8\b\int_{Q_{b_1, b_2}}\lambda\mu\varphi\n^2\eta \n z_s\n z_s\mathrm{dx}\mathrm{ds}\\
\ns\ds\leq -C\b\int_{Q_{b_1, b_2}}\lambda\mu\varphi \th^2\(| \n v_s|^2+ \l^2\mu^2\varphi^2(|v_s|^2+|\n v|^2)+\l^4\mu^4\varphi^4 v^2 \)\mathrm{dx}\mathrm{ds}
\ea\ee

In view of \eqref{0112-ax0}, we only need to estimate $-\b\int_{Q_{b_1, b_2}}\lambda\mu\varphi \th^2| \n v_s|^2\mathrm{dx}\mathrm{ds}$. By using Lemma \ref{lem:5}, we have
\bel{430}\ba{ll}\ds
8\b\int_{Q_{b_1, b_2}}\lambda\mu\varphi\n^2\eta \n z_s\n z_s\mathrm{dx}\mathrm{ds} \\
\ns\ds\leq \frac{C}{-\b\mu^{\frac{1}{2}}}\[||\th\cP v ||^2_{L^2(Q_{b_1, b_2})}+\int_{Q_{b_1,b_2}}\lambda^3\mu^{\frac{7}{2}}\varphi^3\th^2 |v_{s}|^2\mathrm{dx} \mathrm{ds} + \int_{Q_{b_1,b_2}}(A_2+A_3)\mathrm{dx} \mathrm{ds} \]\\
\ns\ds\q+\frac{C}{-\b\mu^{\frac{1}{2}}}\int_{\Sigma_{b_1,b_2}}\(\l^2\mu^3\varphi^2\th^2|\D v|^2+\th^2|\n\D v|^2\)\mathrm{d\Sigma_{b_1,b_2}}  \ea\ee

Plugging \eqref{0112-0v}-\eqref{0112-0d0}, \eqref{1225}, \eqref{i1i548-03} and  \eqref{430} into \eqref{0112-ax0}, we have
\bel{0112-ax01}\ba{ll}\ds
 c\int_{Q_{b_1, b_2}}\th^2\(\lambda^6\mu^8\varphi^6|\n\eta |^8 v^2+\lambda^4\mu^6\varphi^4 |\n\eta |^6 |\n v |^2+\l^3\mu^4\varphi^3|\n\eta|^4|\D v |^2\)\mathrm{dxds}\\
  \ns\ds\q +c\[\int_{Q_{b_1, b_2}}\(\l^3\mu^{\frac{7}{2}}\varphi^3|\n\eta|^4z_s^2+\lambda^3\mu^4\varphi^3|\n\eta|^4|\n\ell\cdot\n z|^2+ \lambda\mu^2\varphi|\n\eta\cdot\n\D z|^2\)\mathrm{dxds} \\
  \ns\ds \q-\int_{Q_{b_1, b_2}}\b\lambda\mu^{\frac{3}{2}}\varphi \th^2|\n\D v|^2\mathrm{dxds}\]\\
 \ns\ds\le C\int_{Q_{b_1, b_2}}\th^2\( |\cP v|^2+\l^4\mu^5\varphi^4|\n v|^2+\l^6\mu^7\varphi^6v^2+\l^2\mu^3\varphi^2|\n^2 v|^2+\b^2\l^3\mu^3\varphi^3v_s^2\)\mathrm{dxds}\].
 \ea\ee

  {\bf Step 4. Removing the $\n \eta$. }
Noting that $z=\th v$, by (\ref{1224-0a}),  we know that there are constants $c_0>0$ and $C>0$ such that
 \bel{0112-bx0}\ba{ll}\ds
  \int_{Q_{b_1, b_2}}\th^2\(\lambda^6\mu^8\varphi^6|\n\eta |^8 v^2+\lambda^4\mu^6\varphi^4 |\n\eta |^6 |\n v |^2+\l^3\mu^4\varphi^3|\n\eta|^4|\D v |^2\)\mathrm{dxds}\\
  \ns\ds\q + \int_{Q_{b_1, b_2}}\(\l^3\mu^{\frac{7}{2}}\varphi^3|\n\eta|^4 z_s^2+\lambda^3\mu^4\varphi^3|\n\eta|^4|\n\ell\cdot\n z|^2\)\mathrm{dxds}\\
  \ns\ds\ge c_0\int_{Q_{b_1, b_2}}\th^2\(\lambda^6\mu^8\varphi^6 v^2+\lambda^4\mu^6\varphi^4  |\n v |^2+\l^3\mu^4\varphi^3\th^2|\D v |^2+\l^3\mu^{\frac{7}{2}}\varphi^3v_s^2\)\mathrm{dxds}\\
  \ns\ds\qq +\int_{Q_{b_1, b_2}}\(\lambda^5\mu^6\varphi^5\th^2|\n\eta\cdot\n v|^2+ \lambda\mu^2\varphi|\n\eta\cdot\n\D v|^2\)\mathrm{dxds}\\
  \ns\ds\qq-C\int_{b_1}^{b_2}\int_{\o_0}\th^2\(\lambda^7\mu^8\varphi^7 v^2+\lambda^5\mu^6\varphi^5  |\n v |^2+\l^3\mu^4\varphi^3(|\D v |^2+v_s^2)\)\mathrm{dxds}.
 \ea\ee

{\bf Step 5. Estimation of lower order terms of $v$ in the R.H.S of \eqref{point-th1-01}.}
  By using the boundary condition $\ds\frac{\pa z}{\pa\nu}=0$ again, we have
 \bel{0112-f0}\ba{ll}\ds
   \int_{Q_{b_1, b_2}}  \l^2\mu^3\varphi^2|\n^2 v|^2\mathrm{dxds}\\
  \ns\ds=-\l^2\mu^3 \int_{Q_{b_1, b_2}} \n^2 v\n(\varphi^2)\cdot\n v\mathrm{dx}\mathrm{ds}-\l^2\mu^3 \int_{Q_{b_1, b_2}}  \varphi^2\n\D v\cdot\n v \mathrm{dxds}\\
  \ns\ds\le \frac{1}{2} \int_{Q_{b_1, b_2}}  \l^2\mu^3\varphi^2|\n^2 v|^2\mathrm{dxds}+C \int_{Q_{b_1, b_2}}  \l^2\mu^3\varphi^2\(\mu^2|\n\eta\cdot\n v|^2+|\D v|^2\) \mathrm{dxds}\\
  \ns\ds +C \int_{Q_{b_1, b_2}}  \l\mu\varphi\(|\n\D v|^2+\l^3\mu^5\varphi^3|\n v|^2\) \mathrm{dxds}
  \ea\ee

 Combining (\ref{0112-ax0})--(\ref{0112-f0}), by Lemma \ref{lem:5},  we conclude that there is $\mu_1>0$ such that for $\mu\ge\mu_1$, there is a $\l_1>0$, for all $\l\ge\l_1$, we have
 \bel{0112-cx02}\ba{ll}\ds
 \int_{Q_{b_1, b_2}}\th^2\(\lambda^6\mu^8\varphi^6 v^2+\lambda^4\mu^6\varphi^4  |\n v |^2+\l^3\mu^4\varphi^3|\D v |^2+ \lambda^3\mu^\frac{7}{2}\varphi^3 v_s^2\)\mathrm{dxds}\\
 \ns\ds\q+\int_{Q_{b_1, b_2}}\th^2\(\lambda\mu^{\frac{3}{2}}\varphi |\n\D v|^2+\lambda\mu^{\frac{3}{2}}\varphi  |\n v_s|^2+\lambda^5\mu^6\varphi^5|\n\eta\cdot\n v|^2+ \lambda\mu^2\varphi|\n\eta\cdot\n\D v|^2\)\mathrm{dxds}\\
 \ns\ds\le C||\th \cP v||^2_{L^2((Q_{b_1, b_2})}\\
 \ns\ds\qq+C \int_{b_1}^{b_2}\int_{\o_0}\th^2\(\lambda^7\mu^8\varphi^7 v^2+\lambda^5\mu^6\varphi^5  |\n v |^2+\l^3\mu^4\varphi^3|\D v |^2+\lambda^3\mu^4\varphi^3v_s^2\)\mathrm{dxds}.
 \ea\ee
Now,  we choose a cut-off function $\chi\in C_0^{\infty} (\o; [0,1])$ so that $\chi=1$ on $\o_0$. Noting that $\cP v-\a v_s=\b v_{ss}+\D^2 v$, then
 \bel{0112-l0}\ba{ll}&\ds
 \l^3\mu^4\varphi^3\chi^2\th^2v(\cP v-\a v_s)\\
 \ns&\ds\ge\frac{\b}{2}\( \l^3\mu^4\varphi^3\chi^2\th^2v_sv\)_s-\frac{\b}{2}\l^3\mu^4\varphi^3\chi^2\th^2v_s^2-\frac{\b}{2}\( \l^3\mu^4\varphi^3\chi^2\th^2\)_s v_sv\\
 \ns&\ds\q+\l^3\mu^4\varphi^3\chi^2\th^2|\D v|^2+\n\cdot(\l^3\mu^4\varphi^3\chi^2\th^2v\n\D v-\l^3\mu^4\varphi^3\chi^2\th^2\n v\D v)\\
 \ns&\ds\q-\n(\l^3\mu^4\varphi^3\chi^2\th^2)\cdot(\n\D v v-\n v\D v).
 \ea\ee
Integrating (\ref{0112-l0}) on $Q_{b_1, b_2}$, noting that $v(b_1)=v(b_2)=0$ and $\ds v=\frac{\pa v}{\pa\nu}=0$ on $\Si_{b_1, b_2}$, we have
 \bel{0112-l1}\ba{ll}\ds
 \int_{Q_{b_1, b_2}}\[-\frac{\b}{2}\l^3\mu^4\varphi^3\chi^2\th^2v_s^2+\l^3\mu^4\varphi^3\chi^2\th^2|\D v|^2\]\mathrm{dxds}\\
 \ns\ds\le C_{\e}\[\int_{Q_{b_1, b_2}}\th^2\(\chi^2|\cP v|^2+\l^7\mu^{8}\varphi^7\chi^2 v^2\)\mathrm{dxds}+\l^5\mu^6\int_{Q_{b_1, b_2}} \th^2\varphi^5\chi^2\th^2|\n v|^2\mathrm{dxds}\]\\
 \ns\ds\q +\e\l\mu^2\int_{Q_{b_1, b_2}}\varphi\th^2\(|\n\eta\cdot\n\D v|^2+\l^2\mu^2\varphi^2\chi^2v_s^2+\l^2\mu^2\varphi^2\chi^2|\D v|^2\)\mathrm{dxds}.
 \ea\ee
  Next, noting that $v=0$ on the boundary, we have
\bel{1225-b02}\ba{ll}\ds
\l^5\mu^6\int_{Q_{b_1, b_2}} \th^2\varphi^5\chi^2\th^2|\n v|^2\mathrm{dxds}\\
 \ns\ds=-\l^5\mu^6\int_{Q_{b_1, b_2}} \th^2\varphi^5\chi^2\th^2\D v v\mathrm{dx}\mathrm{ds}-\l^5\mu^6\int_{Q_{b_1, b_2}}  \n(\th^2\varphi^5\chi^2)\cdot \n v v\mathrm{dx}\mathrm{ds}\\
 \ns\ds\le \frac{1}{2}\l^5\mu^6\int_{Q_{b_1, b_2}}\varphi^5\chi^2\th^2|\n v|^2\mathrm{dxds}+\tilde\e\l^3\mu^4\int_{Q_{b_1, b_2}}\th^2\varphi^3\chi^2|\D v|^2\mathrm{dxds}\\
 \ns\ds\q+C_{\tilde\e}\l^7\mu^{8}\int_{Q_{b_1, b_2}}\th^2\varphi^7\chi^2 v^2\mathrm{dxds}.
\ea\ee
Plugging \eqref{0112-l1} and \eqref{1225-b02} into \eqref{0112-cx02}, taking $\e, \tilde\e>0$ small enough,  we can get the desired result immediately.\endpf

 \ss

\section{Proof of Theorem \ref{thm:3}.}\label{section5-0}

In this section, we are going to prove Theorem \ref{thm:3}.  Since $\beta=0$, $\cP v=\a v_s+\D^2 v$.

 {\bf Proof of Theorem \ref{thm:3}. }
Recall that $z=\th v$ satisfying the assumption $\th(b_1,\cdot)=\th(b_2,\cdot)=0$. Hence, by (\ref{mv-00a0}), we know that
 \bel{0112-0b0}
\int_{Q_{b_1, b_2}} \pa_s M \mathrm{dxds}=0.
 \ee
By means of the boundary conditions $\ds v=\frac{\pa v}{\pa\nu}=0$ and $\b=0$,  proceeding exactly the same analysis as (\ref{1226-h0}) and (\ref{0112-cx02}), we have
 \bel{0112-cx0}\ba{ll}\ds
  \int_{Q_{b_1, b_2}}\th^2\(\lambda^6\mu^8\varphi^6 v^2+\lambda^4\mu^6\varphi^4  |\n v |^2+\l^3\mu^4\varphi^3|\D v |^2\)\mathrm{dxds} \\
 \ns\ds\q+ \int_{Q_{b_1, b_2}}\(\lambda\mu^2\varphi|\n\D z\cdot\n\eta|^2+\lambda^5\mu^6\varphi^5|\n z\cdot\n\eta|^2\)\mathrm{dxds}+\int_{Q_{b_1, b_2}} |\cP_1(z)|^2\mathrm{dxds}\\
  \ns\ds\q+\int_{\Sigma_{b_1,b_2}}\lambda\mu\varphi\th^2\(|\frac{\pa\D v}{\pa\nu}|^2+\l^2\mu^2\varphi^2|\D v|^2 \)\mathrm{dxds}\\
 \ns\ds\le C\[||\th \cP v||^2_{L^2(Q_{b_1, b_2})}+ \int_{Q_{b_1, b_2}} A_4\mathrm{dxds}\]\\
 \ns\ds\qq+C\int_{b_1}^{b_2}\int_{\o_0} \th^2\(\lambda^7\mu^8\varphi^7 v^2+\lambda^5\mu^6\varphi^5  |\n v |^2+\l^3\mu^4\varphi^3|\D v |^2\)\mathrm{dxds}.
 \ea\ee
Here $|\cP_1(z)|^2$ comes from the analysis above Theorem \ref{point-th1}.

In order to get \eqref{i1i2-c10} from \eqref{0112-cx0}, it suffices to estimate
$\ds\int_{Q_{b_1, b_2}}\frac{1}{\l\varphi}\(|\a v_s|^2+|\D^2 v|^2\)\mathrm{dxds}$".

Recalling (\ref{0227-01}) for the definitions of $\cP_j(z)$ ($j=1, 2$), we have that
 \bel{1225-0a0}\ba{ll}\ds
 \int_{Q_{b_1, b_2}}\frac{1}{\l\varphi}|\a z_s|^2\mathrm{dxds}\\
 \ns\ds=\int_{Q_{b_1, b_2}}\frac{1}{\l\varphi}\|\cP_2 (z)+4 \n\ell\cdot\n\D z+4\n^2\ell:\n^2 z+2\D\ell\D z+4 |\n\ell|^2\n\ell\cdot\n z-F z\|^2\mathrm{dxds}\\
 \ns\ds\le  \int_{Q_{b_1, b_2}}\frac{1}{\l\varphi}\|\th \cP v-\cP_1(z)-\cP_r(z)\|^2\mathrm{dxds}\\
 \ns\ds\q+C\[\int_{Q_{b_1, b_2}}\[\l^5\mu^8\varphi^5z^2+\l^5\mu^6\varphi^5|\n\eta\cdot\n z|^2+\l\mu^4\varphi|\n^2 z|^2+\l\mu^2\varphi|\n\eta\cdot\n\D z|^2\]\mathrm{dxds}.
 \ea\ee
 Hence, we have
  \bel{1225-0a1}\ba{ll}\ds
 \int_{Q_{b_1, b_2}}\frac{1}{\l\varphi} \th^2\( |\a v_s|^2+ | \D^2 v|^2\)\mathrm{dxds}\\
 \ns\ds\le \int_{Q_{b_1, b_2}}\frac{C}{\l\varphi}\(|\cP_1 (z)|^2+|\th \cP v|^2\)\mathrm{dxds}\\
 \ns\ds\q + C\int_{Q_{b_1, b_2}}\[\l^5\mu^6\varphi^5|\n\eta\cdot\n z|^2+\l\mu^2\varphi|\n\eta\cdot\n\D z|^2+A_1+A_3\]\mathrm{dxds}.\\
 \ea\ee
 Next, we compute
 \bel{1219-f2}\ba{ll}&\ds
 \l\mu^2\varphi\th^2|\n\D v|^2\\
 \ns&\ds= \n\cdot\(\l\mu^2\varphi\th^2\D v\n\D v\)-\n (\l\mu^2\varphi\th^2)\cdot\n\D v\D v-\l\mu^2\varphi\th^2\D v\D^2 v\\
 \ns&\ds\le \n\cdot\(\l\mu^2\varphi\th^2\D v\n\D v\)+\frac{1}{2} \l\mu^2\varphi\th^2|\n\D v|^2+C\l^3\mu^4\varphi^3\th^2|\D v|^2+\frac{1}{\l\varphi}\th^2|\D^2 v|.
 \ea\ee
 Hence, integrating \eqref{1219-f2} on $Q_{b_1, b_2}$, we have
  \bel{1225-0aj}\ba{ll}\ds
  \int_{Q_{b_1, b_2}} \l\mu^2\varphi\th^2|\n\D v|^2\mathrm{dxds}&\ds\le C\int_{\Sigma_{b_1,b_2}}\mu\th^2\(|\frac{\pa\D v}{\pa\nu}|^2+\l^2\mu^2\varphi^2|\D v|^2 \)\mathrm{dxds}\\
  \ns&\ds\qq+C \int_{Q_{b_1, b_2}} \(\l^3\mu^4\varphi^3\th^2|\D v|^2+\frac{1}{\l\varphi}\th^2|\D^2 v|\)\mathrm{dxds}.
  \ea\ee

Combining (\ref{0112-cx0}), (\ref{1225-0a1}) and  (\ref{1225-0aj}), noting that  $z=\th v$, we conclude that there is $\mu_1>0$ such that for $\mu\ge\mu_1$, there is a $\l_1>0$, for all $\l\ge\l_1$, we have
 \bel{0112-x1}\ba{ll}\ds
   \int_{Q_{b_1, b_2}} \th^2\[\l^6\mu^8\varphi^6v^2+\l^4\mu^6\varphi^4|\n v|^2+\l^3\mu^4\varphi^3|\D v|^2+\l\mu^2\varphi|\n\D v|^2\]\mathrm{dxds}\\
  \ns\ds\le C\[ \int_{Q_{b_1, b_2}} \th^2|\cP v|^2\mathrm{dxds}+\int_{b_1}^{b_2}\int_{\o_0}\th^2\[\l^7\mu^8\varphi^7v^2+\l^5\mu^6\varphi^5|\n v|^2+\l^3\mu^4\varphi^3|\D v|^2\]\mathrm{dxds}.
  \ea\ee

Finally, proceeding the same analysis as (\ref{0112-l0})--(\ref{1225-b02}), we can prove that
 $$ \int_{b_1}^{b_2}\int_{\o_0}\th^2\[\l^5\mu^6\varphi^5|\n v|^2+\l^3\mu^4\varphi^3|\D v|^2\]\mathrm{dxds}\le C\l^7\mu^8\int_{b_1}^{b_2}\int_{\o}\varphi^7\th^2 v^2\mathrm{dxds},$$
 which yields  the desired result immediately.\endpf

\section{Proof of Theorem \ref{0tt2}}\label{section6-0}

Before giving the proof of Theorem \ref{0tt2}, we fist establish the following interpolation inequality, which will play a key role in the proof of resolvent estimate.
 \bt\label{1202-t1}
Let $d$ satisfying (\ref{cond1.1}). Then there are constants $\l^*>0$ and
$C=C(\l^*)>0$ such that for all $\l\geq\l^*$ and
any solution $q$ satisfying
 \bel{1202-s1}\left\{\ba{ll}\ds
 q_{ss}-\D^2 q+id q_s
=q^0 & \hb { in } (-2, 2)\t \O,\\
\ns\ds q=\frac{\pa q}{\pa\nu}=0,\q \mbox{ or }\q q=\D q=0& \hb { on }(-2, 2)\t \pa\O.
 \ea\right.\ee
with $q^0\in L^2((-2,2)\t\O)$ , it holds that
$$\ba{ll}\ds |q|_{H^1(-1,1; L^2(\O))}+|q|_{L^2(-1,1; H^2(\O))}\\
  \ns\ds\le
 Ce^{C\l}\big(|\!|q^0|\!|_{L^2((-2,2)\t\O)}+|\!|q|\!|_{L^2((-2,2)\t\o)}+
 |\!|d q_s|\!|_{L^2((-2,2)\t \O))} \big)\\
 \ns\ds\q+Ce^{-\l}|\!|q|\!|_{ {H^1(-1,1; L^2(\O))}+|q|_{L^2(-1,1; H^2(\O))}}.
  \ea$$
 \et

 \ms

{\it Proof of Theorem \ref{1202-t1}. } We divide the proof into four steps.

\ss

{\it Step 1. }  Noting that there is no boundary conditions
for $q$ at $s=\pm2$ in the system
(\ref{1202-s1}), thus we need to introduce a
cut-off function $\eta=\eta(s)\in
C_0^\infty(-b,b)\subset C_0^\infty(\dbR)$ such
that
 \bel{1202-as7}\left\{\ba{ll}\ds
 0\le\eta(s)\le1,  &|s|<b,\\
 \ns\ds \eta(s)=1,&|s|\le b_0.
  \ea\right.\ee
Here $1<b_0<b\le2$ are given as follows:
\bel{1202-as2} b\=\sqrt{1+\frac{1}{\mu}\ln
(2+e^\mu)},\qq b_0\=\sqrt{b^2-{\frac{1}{\mu}}\ln
\left(\frac{1+e^\mu}{e^\mu}\right)}, \ee
where $\mu$ is the parameter  which is large enough such that
$b<2$. For parameter $\l,\mu>1$, we define the weight functions $\th$ and $\phi$ as follow:
 \bel{cas4}
 \th=e^{\ell}, \q \ell=\l\phi,\q \phi=e^{\mu\psi},\q \psi(s,x)={\eta(x)\over
 |\eta|_{L^\infty(\O)}}+b^2-s^2
  \ee
where $\eta\in
C^2(\oO)$ satisfying (\ref{1224-0a}).

Put
 \bel{1202-as8}
\hat q=\eta q.
 \ee
Noting that $\varphi$ does not depend on $x$, it
follows from (\ref{1202-s1}) that
 \bel{1202-as9}\left\{\ba{ll}\ds
\hat q_{ss}-\D^2 \hat q=\eta_{ss}q+2\eta_sq_s+\eta(q^0-i d q_s)& \mbox{ in } (-2,2)\t\O, \\
\ns\ds\hat q=\frac{\pa \hat q}{\pa\nu}=0, \ \mbox{ or }\  \hat q=\D \hat q=0&\mbox{ on
}(-2,2)\t\pa\O.
 \ea\right.
 \ee

In Theorem \ref{thm:2},  by taking $\b=-1$, $[b_1, b_2]=[-b, b]$, $v$ replaced by $\hat q$, we conclude that there exists
constant  $\mu_0>0$, such that for all $\mu\ge\mu_0$, there exist constants $C=C(\mu)>0$ and $\l_0(\mu)>0$, for all $\l\ge\l_0$, it
holds that
  \bel{0104-c1}\ba{ll}\ds
 \int_{-2}^{2}\int_\O\th^2 \(\lambda^6\mu^8\varphi^6  \hat q^2+\lambda^4\mu^6\varphi^4 |\n \hat q |^2
+\l^3\mu^3\varphi^3 (|\D \hat q |^2+|\hat q_s|^2)\)\mathrm{dxds}\\
 \ns\ds\le C\[\int_{-2}^2\int_\O\th^2\|\eta_{ss}q+2\eta_sq_s+\eta(q^0-i d q_s)\|^2 \mathrm{dxds}+\lambda^7\mu^8\int_{-2}^{2}\int_\o\varphi^7\th^2 \hat q^2\mathrm{dxds} \].
 \ea\ee
Recalling (\ref{cas4}) for the definition of $\phi$, it is easy to see that
 \bel{0107-as5}\left\{\ba{ll}\ds
 \phi(s,\cd)\ge 2+e^\mu , &\hb{ for any $s$ satisfying }|s|\le 1,\\
 \ns\ds \phi(s,\cd)\le 1+e^{\mu}, &\hb{ for any $s$ satisfying }b_0\le |s|\le
 b.
 \ea\right.\ee
Denoting $c_0=2+e^{\mu}>1$, and note that $b_0\in (1,b)$. By using
(\ref{0104-c1}) and (\ref{0107-as5}), note that $\hat q=\varphi q$, we get
 \bel{as217}\ba{ll}\ds
\l e^{2\l c_0}\int_{-1}^1\int_{\O}(|\n
 q|^2+|q_s|^2+|\D q|^2+|q|^2)\mathrm{dxds}\\[3mm]
\ns\ds\le
Ce^{C\l}\left\{\int_{-2}^2\int_{\O}(|q^0|^2+|d q_s|^2)\mathrm{dxds}+\int_{-2}^2\int_{\o}q^2\mathrm{dxds}\right\}\\[3mm]
\ns\ds\q+Ce^{2\l(c_0-1)}\int_{(-b,-b_0)\bigcup (b_0,
b)}\int_{\O}(|q|^2+|q_s|^2)\mathrm{dxds}.
 \ea\ee
This completes the proof of Theorem \ref{1202-t1}. \endpf

\ms

{\it Proof of Theorem \ref{0tt2}. } Fix $F=(g^0,g^1)\in \cH$
and $U_0 =(z^0,z^1)\in D(\cA)$. Then,
 \bel{6a1}
(\cA-i\gamma I)U_0=F
 \ee
is equivalent to
 \bel{6a2}\left\{\ba{ll}\ds
 -i\gamma z^0+z^1=g^0  & \hb { in } \O,\\
 \ns\ds
-\D^2 z^0-[d +i\gamma]z^1=g^1   & \hb {
in }
 \O,\\
\ns\ds z^0=\frac{\pa z^0}{\pa \nu}=0 \ \mbox{ or } z^0=\D z^0=0&\hb{ on } \pa\O.\\
 \ea\right.\ee
 By (\ref{6a2}), we conclude that $z^0$ satisfying
 \bel{6a3}
-\D^2z^0+\gamma^2 z^0-i\omega d z^0=[d +i\gamma] g^0+g^1  \mbox { in } \O.
\ee
 For any $s\in\dbR$, set  $q=e^{-\gamma s}z^0$. Then $q$  solves the following equation:
 \bel{6a4}\left\{\ba{ll}\ds
 q_{ss}-\D^2 q+id q_s
=\big[i \gamma g^0+d g^0+g^1\big]e^{-\gamma s} & \hb { in } \dbR\t \O,\\
\ns\ds q=\frac{\pa q}{\pa\nu}=0, \ \mbox{ or } q=\D q=0& \hb { on }\dbR\t \pa\O.
 \ea\right.\ee
Noting that $q=e^{-\gamma s}z^0$ and $|\gamma|>1$, we have
\bel{6a7-0}\ba{ll}\ds
|\!|z^0|\!|_{H^2(\O)}\le
Ce^{C|\gamma|}\(||q||_{H^1(-1,1; L^2(\O))}+||q||_{L^2(-1,1; H^2(\O))}\).
 \ea\ee
Similarly, it is easy to get that
\bel{6a7}
 ||q||_{H^1(-2,2; L^2(\O))}+||q||_{L^2(-2,2; H^2(\O))}
 \le C(|\gamma|+1)e^{C|\gamma|}||z^0||_{H^2(\O)},\ee
 \bel{6a43}
 ||q||_{L^2((-2,2)\t\o)}+||dq_s||_{L^2(-2,2; H^2(\O))}\le C(1+|\gamma|)e^{C|\gamma|}\(||z^0||_{L^2(\o)}+||dz^0||_{L^2(\O)}\).
  \ee
 Applying Theorem \ref{1202-t1} to the equation (\ref{6a4}), by
(\ref{6a7}), we obtain that
 \bel{06a5}\ba{ll}\ds
||z^0||_{H^2(\O)}\le
Ce^{C|\gamma|}\big(||g^0||_{H^2(\O)}+||g^1||_{L^2(\O)}+||z^0||_{L^2(\o)}+||dz^0||_{L^2(\O)}\big).
 \ea\ee
Multiplying the first equation of (\ref{6a3})
by $2\overline z^0$ and integrating it on $\O$,
we find that
\bel{n6a8}\ba{ll}\ds \int_{\O}2\big[\D^2
z^0+\l^2z^0+i\gamma d z^0\big] \overline
z^0\mathrm{dx}=2\l^2\int_\O|z^0|^2\mathrm{dx}+2\int_\O |\D
z^0|^2\mathrm{dx}+2i\gamma \int_{\O}d |z^0|^2\mathrm{dx}.
 \ea\ee
Taking the imaginary part in  (\ref{n6a8}), by (\ref{6a4}),
we have
 \bel{6a9}\ba{ll}\ds
2|\gamma|\int_{\O}d |z^0|^2\mathrm{dx} \le C||i\gamma
g^0+d g^0+g^1||_{L^2(\O)}||z^0||_{L^2(\O)}.
 \ea\ee
Combining (\ref{06a5}) and (\ref{6a9}), and
noting that $d \ge d_0>0$ in $\o$, we conclude that there exists a constant $C>0$ such that
 \bel{6a10}\ba{ll}\ds
||z^0||_{H^2(\O)}\le
Ce^{C|\gamma|}\[||g^0||_{H^2(\O)}+||g^1||_{L^2(\O)}\].
 \ea\ee

Recalling that $z^1=g^0+i\gamma z^0$, it follows
  \bel{6a12}\ba{ll}\ds
||z^1||_{L^2(\O)}\le
Ce^{C|\gamma|}\big(||g^0||_{H^2(\O)}+||g^1||_{L^2(\O)}\big).
 \ea\ee
Finally, by (\ref{6a10})--(\ref{6a12}), we finish the proof of Theorem \ref{0tt2}.\endpf

\ms
\appendix
\section{Appendix: Proofs of Lemma \ref{lem:1} and Lemma \ref{lem:02}}

\subsection{Proof of Lemma \ref{lem:1}}
{\bf Proof of Lemma \ref{lem:1}. }
 For  convenience, we recalled
 \bel{0624-0a}\left\{\ba{ll}\ds
 \vec{E}=4\n A+4\D\ell\n\ell,\q F=\Delta^2 \ell-2\n\cdot(A\n\ell),\\
 \cP_3(z)=\Phi z+ Z\n^2 z:\n^2\eta,\\
 J_1(z)=\D^2 z+A^2 z+2A\D z+4\n^2 z\n\ell\n\ell+\vec{E}\cdot\n z,\\
 \ns\ds J_2(z)=-4 \n\ell\cdot\n\D z-4\n^2\ell:\n^2 z-2\D\ell\D z-4 A\n\ell\cdot\n z+Fz.
 \ea\right.\ee
 We divide the proof into three steps.

\ms

  {\bf Step 1. Estimation of $\ds\a z_s\(z_{ss}+J_1(z)\)$. }  First, by elementary calculus, we have
 \bel{i11i21}\ba{ll}&\ds
 \a\b z_s z_{ss}+\a z_s J_1(z)\\
 \ns&\ds=\frac{1}{2}\(\a\b z_s^2\)_s+\n\cdot\(\a z_s\n\D z-\a\n z_s\D z\)+\frac{1}{2}\(\a|\D z|^2+\a A ^2 z^2\)_s-\a A A_s z^2\\
\ns&\ds\q+\n\cdot(2\a A z_s\n z)-\(A\a |\n z|^2\)_s+\a A_s|\n z|^2-2\a z_s\n A\cdot\n z\\
 \ns&\ds\q +\n\cdot\(4\a z_s(\n\ell\cdot\n z)\n\ell\)-4\a(\n\ell\cdot\n z)(\n\ell\cdot\n z_{s})\\
 \ns&\ds\q-4\a z_s(\n^2\ell\n\ell)\cdot\n z-4\a\D\ell z_s(\n\ell\cdot\n z)+\a  z_s\vec{E}\cdot\n z.
 \ea\ee

  Remember that $A=|\n\ell|^2$, then
  \bel{na}
  \n A=2\n^2\ell\n\ell.
  \ee
 Therefore, by (\ref{na}) and (\ref{0227-01}), it is easy to see that
   \bel{1208-c}\ba{ll}\ds
   -2\a z_s\n A\cdot\n z-4\a z_s(\n^2\ell\n\ell)\cdot\n z-4\a\D\ell z_s(\n\ell\cdot\n z)=-\a  z_s\vec{E}\cdot\n z.
   \ea\ee
  Noting that
  \bel{na1}\ba{ll}\ds
  -4\a(\n\ell\cdot\n z)(\n\ell\cdot\n z_{s})=-2\a\(|\n\ell\cdot\n z|^2\)_s+4\a(\n\ell\cdot\n z)(\n\ell_s\cdot\n z).
  \ea\ee
Then, combining (\ref{i11i21}), (\ref{1208-c}) and (\ref{na1}), by (\ref{1211-0a0}), (\ref{0421-02}) and (\ref{A11-12-13}),  we have
\bel{i1i21-0}\ba{ll}\ds
 \a\b z_s z_{ss}+\a z_sJ_1(z)&\ds=\pa_s M_1+\n\cdot V_1-\a A A_s z^2+\a A_s|\n z|^2+4\a(\n\ell\cdot\n z)(\n\ell_s\cdot\n z)\\
 \ns&\ds\ge \pa_s M_1+\n\cdot V_1-CA_1.
   \ea\ee

\ms

{\bf Step 2. } Let us estimate ``$\ds \b z_{ss}J_2(z)$". Firstly, we
 notice that
  \bel{1208-2a}\ba{ll}\ds
  -4 \b z_{ss}(\n\ell\cdot\n\D z)=-4\b\n\cdot\( z_{ss} \n\ell\D z\) +4\b\n z_{ss}\cdot \n\ell\D z+4\b\D\ell\D z z_{ss}.
  \ea\ee
Moreover, we have
  \bel{ddt11}\ba{ll}\ds
 4\b\n z_{ss}\cdot \n\ell\D z\\
  \ns\ds=4\(\b\n\ell\cdot\n z_s\D z\)_s-4\b\n\cdot\(\n\ell_s\D z z_s\)+4\b\D\ell_s\D z z_s+4\b\n\ell_s\cdot\n\D z z_s\\
  \ns\ds\q+\n\cdot\[-4\b\(\n\ell\cdot\n z_s\)\n z_s+2\b\n\ell|\n z_s|^2\]+4\b\n^2\ell\n z_s\n z_s-2\b\D\ell|\n z_s|^2
  \ea\ee
and
   \bel{dr24}\ba{ll}\ds
 4\b\D\ell\D z z_{ss}&\ds=4\(\b\D\ell\D z z_s\)_s-4\b\D\ell_s\D z z_s-4\b\n\cdot\(\D\ell\n z_s z_s\)\\
 \ns&\ds\q+4\b\D\ell|\n z_s|^2+2\b\n\cdot\(\n\D\ell| z_s|^2\)-2\b\D^2\ell| z_s|^2
  \ea\ee
Thus, by (\ref{1208-2a})--(\ref{dr24}), we have
 \bel{1208-2b}\ba{ll}\ds
  -4 \b z_{ss}(\n\ell\cdot\n\D z)\\
  \ns\ds=\n\cdot\[ -4\b z_{ss} \n\ell\D z-4\b\n\ell_s\D z z_s-4\b\(\n\ell\cdot\n z_s\)\n z_s+2\b\n\ell|\n z_s|^2\\
  \ns\ds\q-4\b\D\ell\n z_s z_s+2\b\n\D\ell| z_s|^2\]+\(4\b\n\ell\cdot\n z_s\D z+4\b\D\ell\D z z_s\)_s\\
  \ns\ds\q+4\b\n\ell_s\cdot\n\D z z_s+4\b\n^2\ell\n z_s\n z_s+2\b\D\ell|\n z_s|^2-2\b\D^2\ell z_s^2.
  \ea\ee
 Next,
 \bel{1208-b}\ba{ll}
 -4\b z_{ss}\n^2\ell:\n^2 z&\ds=-4\b\n\cdot\(\n^2\ell\n z z_{ss}\)+4\b\n\D\ell\cdot\n z z_{ss}+4\b\n^2\ell \n z_{ss}\n z\\
\ns&\ds=-4\b\n\cdot\(\n^2\ell\n z z_{ss}\)+4\b\(\n\D\ell\cdot\n z z_s\)_s-4\b\n\D\ell_s\cdot\n z z_s\\
  \ns&\ds\q-2\b\n\cdot\(\n\D\ell| z_s|^2\)+2\b\D^2\ell| z_s|^2+\(4\b\n^2\ell\n z\n z_s\)_s\\
   \ns&\ds\q-4\b\n^2\ell\n z_s\n z_s-2\b\(\n^2\ell_s\n z\n z\)_s+2\b\n^2\ell_{ss}\n z\n z
 \ea\ee
and
\bel{1208-3a}\ba{ll}\ds
 -2\b z_{ss}\D\ell\D z&\ds= -2\b \(z_s\D\ell\D z\)_s+2\b\n\cdot\(\D\ell_s\n z z_s\) -2\b z_s\n\D\ell_s\cdot\n z\\
 \ns&\ds\q-\b\(\D\ell_s|\n z|^2\)_s+\b\D\ell_{ss}|\n z|^2+2\b\n\cdot\(z_s\D\ell\n z_s \)\\
 \ns&\ds\q-2\b\D\ell|\n z_s|^2-\b \n\cdot\(\n\D\ell z_s^2\)+\b\D^2\ell z_s^2.
 \ea\ee

Further,
 \bel{1208-4a}\ba{ll}&\ds
 \b z_{ss}\( -4 A\n\ell\cdot\n z\)\\
\ns&\ds= -4\b \(A\n\ell\cdot\n z z_s\)_s+4\b (A\n\ell)_s\cdot\n z z_s+2\b \n\cdot \(A\n\ell z_s^2\)-2\b\n\cdot \(A\n\ell\)z_s^2.
 \ea\ee
and
\bel{1208-5a}\ba{ll}&\ds
\b F z_{ss} z=\(\b F z_s z\)_s-\frac{1}{2}\b\(F_s z^2\)_s+\frac{1}{2}\b F_{ss} z^2-\b F z_s^2.
 \ea\ee
  Together with (\ref{1208-2b})--(\ref{1208-5a}), by (\ref{1211-0a0}), (\ref{0421-02}) and (\ref{A11-12-13}), we have
  \bel{1208-06a}\ba{ll}\ds
 \b z_{ss}J_2(z)&\ds=\n\cdot V_2+\pa_s M_2+4\b\n\ell_s\cdot\n\D z z_s-6\b\n\D\ell_s\cdot\n z z_s\\
   \ns&\ds\q+2\b\n^2\ell_{ss}\n z\n z+\b\D\ell_{ss}|\n z|^2+\b\[\D^2\ell -2\n\cdot(A\n\ell)\]z_s^2\\
   \ns&\ds\q+4\b (A\n\ell)_s\cdot\n z z_s+\frac{1}{2}\b F_{ss} z^2-\b F z_s^2\\
   \ns&\ds\ge  \n\cdot V_2+\pa_s M_2-C\(A_1+A_2+A_3+\l\mu\varphi|\n\eta\cdot\n\D z|^2\).\\
 \ea\ee

{\bf Step 3. } Finally, let us estimate ``$z_{ss}\cP_3(z)$".
In view of condition \eqref{cd1}, $\eta$ is independent with $s$, then  we have
 \bel{1208-8a}\ba{ll}&\ds
 \b z_{ss}\cP_3(z)=\b z_{ss}\(\Phi z+Z\n^2\eta: \n^2 z\)\\
 \ns&\ds=\(\b\Phi z_s z\)_s-\frac{1}{2}\b\(\Phi_s z^2\)_s+\frac{1}{2}\b \Phi_{ss} z^2-\b\Phi z_s^2+\b \(Z z_{s}\n^2 \eta:\n^2 z\)_s\\
 \ns&\ds\q-\b Z_s z_{s}\n^2 \eta:\n^2 z-\b Z z_{s}\n^2 \eta:\n^2 z_s.
\ea\ee

On the one hand,
 \bel{1208-9a}\ba{ll}\ds
  -\b Z_s z_{s}\n^2 \eta:\n^2 z\\
  \ns\ds= -\b \n\cdot\(Z_s z_{s}\n^2 \eta \n z\)+\b z_{s}\n^2 \eta\n Z_s\n z+\frac{1}{2}\b \(Z_s\n^2 \eta\n z\n z\)_s\\
 \ns\ds\q-\frac{1}{2}\b Z_{ss}\n^2 \eta\n z\n z+\b Z_s z_{s} \(\n\D \eta\cdot\n z\)
 \ea\ee
On the other hand,
\bel{1208-10a}\ba{ll}\ds
 -\b Z z_{s}\n^2 \eta:\n^2 z_s\\
 \ns\ds= -\b \n\cdot\(Z z_{s}\n^2  \eta\n z_s\)+\b Z\n^2 \eta\n z_{s}\n z_s+\frac{1}{2}\b \n\cdot\((\n^2 \eta\n Z) z_s^2\)\\
 \ns\ds\q-\frac{1}{2}\b \n\cdot\((\n^2 \eta\n Z)\)z_s^2+\frac{1}{2}\b \n\cdot\(Z \n\D \eta z_s^2\)-\frac{1}{2}\b \n\cdot\(Z \n\D \eta\) z_s^2.
 \ea\ee
Combining (\ref{1208-8a}) and  (\ref{1208-9a})--(\ref{1208-10a}), we have
  \bel{1208-11a}\ba{ll}\ds
 \b z_{ss}\cP_3(z)\ds&\ds=\pa_sM_3+\n\cdot V_3-\frac{1}{2}\b\[2\Phi + \n\cdot\((\n^2 \eta\cdot\n Z)\)+ \n\cdot\(Z \n\D \eta\) \]z_s^2\\
 \ns&\ds\q+\frac{1}{2}\b \Phi_{ss} z^2+\b z_{s}\n^2 \eta\n Z_s\n z-\frac{1}{2}\b Z_{ss}\n^2 \eta\n z\n z\\
 \ns&\ds\q+\b Z_s z_{s} \(\n\D \eta\cdot\n z\)+\b Z\n^2 \eta\n z_{s}\n z_s\\
\ea\ee
Recalling   (\ref{1211-0a0}), (\ref{0421-02}) and (\ref{A11-12-13}),  we obtain \eqref{i1i21-p1p2}, which concludes  the proof. \endpf.

\subsection{Proof of Lemma \ref{lem:02}}
{\bf Proof of Lemma \ref{lem:02}. } The proof is long, we divide our proof into several steps.

{\bf Step 1. }Let us estimate ``$\ds \D^2 z J_2(z)$".
 By (\ref{0624-0a}), we know that
  \bel{1209-0a}\ba{ll}\ds
  \D^2 z J_2(z)=\D^2 z\[-4 \n\ell\cdot\n\D z-4\n^2\ell:\n^2 z-2\D\ell\D z-4 A\n\ell\cdot\n z+Fz\].
  \ea\ee
First,
 \bel{1209-01}\ba{ll}\ds
 -4\D^2 z\n\ell\cdot\n\D z\\
 \ns\ds=\n\cdot\(-4\n\D z(\n\ell\cdot\n\D z)+2\n\ell|\n\D z|^2\)+4\n^2\ell\n\D z\n\D z-2\D\ell|\n\D z|^2.
 \ea\ee
Next,
  \bel{1209-02}\ba{ll}\ds
 -4\D^2 z\n^2\ell:\n^2 z\\
 \ns\ds= -4\n\cdot\(\n\D z(\n^2\ell:\n^2 z)\)+ 4\n(\n^2\ell:\n^2 z)\cdot\n\D z\\
 \ns\ds=-4\n\cdot\(\n\D z(\n^2\ell:\n^2 z)\)+4\sum_{i, j, k=1}^n\(\ell_{x_ix_j} z_{x_jx_k} \D z_{x_k}\)_{x_i}\\
  \ns\ds\q-4\n \cdot\(\n^2 z \n\D\ell\D z\)+4\n\D z\cdot \n\D\ell \D z+4\n^2 z:\n^2 \D\ell \D z\\
  \ns\ds\q-4\sum_{i, j, k=1}^n\(\ell_{x_ix_j} z_{x_jx_k} \D z_{x_i}\)_{x_k}+4\n^2\ell\n\D z\n\D z+8\sum_{i, j, k=1}^n\ell_{x_ix_jx_k} z_{x_jx_k}\D z_{x_i}.
  \ea\ee

Further,
 \bel{1209-030}\ba{ll}\ds
 -2\D^2 z\D\ell\D z=  -\n\cdot\(2\D\ell\D z\n\D z\)+2\n\D\ell\cdot\n\D z\D z+2\D\ell|\n\D z|^2.
 \ea\ee
and
 \bel{1209-03}\ba{ll}\ds
 -4 A\D^2 z \(\n\ell\cdot\n z\)=-4\n\cdot\(\n\D z(A\n\ell\cdot\n z)\)+4\n\D z\cdot\n (A\n\ell\cdot\n z).
 \ea\ee
We observe that

 \bel{1209-04}\ba{ll}\ds
 4\n\D z\cdot\n (A\n\ell\cdot\n z)\\
 \ns\ds=4\n\D z\cdot\n A(\n\ell\cdot\n z)+4A\n^2\ell\n\D z \n z+4A\n^2 z\n\D z\n\ell\\
 \ns\ds=4\n\cdot\[\D z\n A(\n\ell\cdot\n z)\]-4\D A(\n\ell\cdot\n z)\D z-4\D z\n^2 z\n A\n\ell\\
 \ns\ds\q-4\D z\n^2\ell\n A\n z + 4\n\cdot\( A\n^2\ell\n z\D z\)\\
\ns\ds\q-4A \n^2\ell:\n^2 z\D z-4{A(\n\D\ell )}\cdot\n z\D z\\
\ns\ds\q-4\n^2\ell\n A\n  z\D z+4\n\cdot\[A\n^2 z \(\n^2 z\n\ell\)\]-4A\n^2\ell\n^2 z \n^2 z.\\
\ns\ds\q-4\(\n^2 z\n\ell\)\( \n^2 z\n A\)-2\n\cdot\(A\n\ell|\n^2 z|^2\)+2\n\(A\n\ell\)|\n^2 z|^2
\ea\ee

A short calculation shows that
 \bel{djh2}\ba{ll}&\ds
-4A\n^2\ell\n^2 z\n^2 z\\
\ns&\ds=-4\n\cdot \(A\n^2\ell(\n z\cdot\n^2 z)\)+2\n\cdot\((\n A\cdot\n^2\ell)|\n z|^2\)\\
\ns&\ds\q-4\n\cdot(\n A\n^2\ell)|\n z|^2+4\n\cdot \(Az \n^2z \n\D\ell\)-4z\n^2 z\n A\n\D\ell \\
\ns&\ds\q-4Az\n\D z\cdot \n\D\ell-4Az\n^2 z:\n^2\D\ell +4\n\cdot\(A\n^2\ell:\n^2 z\n z\)\\
\ns&\ds\q-4(\n A\cdot\n z)(\n^2\ell:\n^2 z)-4A\sum_{i, j, k=1}^n\ell_{x_ix_jx_k} z_{x_jx_k}z_{x_i}-4A\n^2\ell:\n^2 z\D z
 \ea\ee

and
\bel{djf2}\ba{ll}&\ds
-4\D z\n^2 z\n A\n\ell\\
\ns&\ds=-4\n\cdot\(\D z\n z\cdot\n A\n\ell\)+4\D\ell\D z\n z\cdot\n A+4\D z\n^2 A\n z\n\ell\\
\ns&\ds\q+4\n\cdot\(\n^2 z\n\ell\n z\cdot\n A\)-4\n^2 z:\n^2\ell\n z\cdot\n A-4\n^2 z\n\ell\(\n^2 A\n z\)\\
\ns&\ds\q -4(\n^2z\n\ell)\cdot(\n^2z \n A)
\ea\ee

Next,
\bel{dj55}\ba{ll}&\ds
Fz\D^2 z=\n\cdot\(F \n\D z z-F\D z\n z\)+F |\D z |^2-\n F \cdot\n\D z z+\n F\cdot\n z \D z\\
\ea\ee

 Combining (\ref{1209-0a}), (\ref{1209-01})--(\ref{1209-04}), we end up with
  \bel{1209-06}\ba{ll}\ds
  \D^2 z J_2(z)&\ds=\n\cdot V_4+8\n^2\ell\n\D z\n\D z+2\n\cdot( A\n\ell)|\n^2 z |^2+F |\D z |^2\\
  \ns&\ds\q-8\(\n^2 z\n\ell\)\( \n^2 z\n A\)-8 A\D z(\n^2\ell:\n^2 z)+R_1,
  \ea\ee
where
 \bel{1210-R1}\ba{ll}\ds
 R_1=6\n\D z\cdot \n\D\ell \D z+4\n^2 z:\n^2 \D\ell \D z+8\sum_{i, j, k=1}^n\ell_{x_ix_jx_k} z_{x_jx_k}\D z_{x_i}-4\D A(\n\ell\cdot\n z)\D z\\
\ns\ds\qq\q-4{A(\n\D\ell )}\cdot\n z\D z-8\n^2\ell\n A\n  z\D z\\
\ns\ds\qq\q-4z\n^2 z\n A\n\D\ell-4Az\n\D z\cdot \n\D\ell-4Az\n^2 z:\n^2\D\ell \\
\ns\ds\qq\q-4A\sum_{i, j, k=1}^n\ell_{x_ix_jx_k} z_{x_jx_k}z_{x_i}+4\D\ell\D z\n z\cdot\n A\\
\ns\ds\qq\q+4\D z\n^2 A\n z\n\ell-8\n^2 z:\n^2\ell\n z\cdot\n A-4\n^2 z\n\ell\(\n^2 A\n z\)\\
\ns\ds\qq\q-\n F \cdot\n\D z z+\n F\cdot\n z \D z.
 \ea\ee

\ms

{\bf  Step 2. Estimation of $A^2 z J_2(z)+2A\D zJ_2(z)$. }

\ms

Noting that
\bel{fg}\ba{ll}
\n\cdot(\n^2 z\n\ell)=\n\ell\cdot \n\D z+\n^2\ell:\n^2 z,
\ea\ee
then by (\ref{0624-0a}), we know that
  \bel{i11-0}\ba{ll}\ds
A^2  z J_2(z)&\ds =A^2 z\(-4 \n\ell\cdot\n\D z-4\n^2\ell:\n^2 z-2\D\ell\D z-4 A\n\ell\cdot\n z+Fz\)\\
 \ns&\ds=-4 A^2 z\n\cdot(\n^2 z\n\ell)-2A^2\D\ell z\D z-4A^3 z\n\ell\cdot\n z+A^2Fz^2\\
 \ns&\ds=\n\cdot\[-4 A^2 z\n^2 z\n\ell+2A^2\n\ell |\n z|^2-2A^2\D\ell z\n z\]\\
 \ns&\ds\q+4\n\cdot\[ z\n A^2\(\n z\cdot \n\ell\)\]-4 z\(\n z\cdot\n \ell\) \D A^2-4 z\(\n^2 \ell\n z\n A^2\)\\
 \ns&\ds\q -4\(\n z\cdot\n \ell\)\(\n z\cdot\n A^2\)-2\n A^2\cdot\n\ell |\n z|^2+2\n(A^2\D\ell)\cdot\n z z\\
 \ns&\ds\q-2\n\cdot\(A^3 \n\ell z^2\)+2\n\cdot\(A^3\n\ell\) z^2+A^2Fz^2.
 \ea\ee
Similarly,

 \bel{i13i22}\ba{ll}&\ds
2A\D z J_2(z)\\
\ns&\ds =2A\D z\(-4 \n\ell\cdot\n\D z-4\n^2\ell:\n^2 z-2\D\ell\D z-4A\n\ell\cdot\n z+Fz\)\\
\ns&\ds=-4\n\(A\n\ell|\D z|^2\)+4\n\(A\n\ell\)|\D z|^2-8A\D z \n^2\ell:\n^2 z\\
 \ns&\ds-4A\D\ell|\D z|^2-8\n\cdot\(A^2\n z(\n\ell\cdot\n z)\)+8(\n A^2\cdot\n z)(\n\ell\cdot\n z)\\
 \ns&\ds\q+8A^2\n^2 \ell\n z\n z+4\n\cdot\(A^2\n\ell |\n z|^2\)-4\n\cdot\( A^2\n\ell\)|\n z|^2\\
 \ns&\ds\q+2\n\cdot\( A F z\n z\)-2A F |\n z|^2-2\n\(A F\)\cdot \n z z.
 \ea\ee

 Combining (\ref{i11-0}) and (\ref{i13i22}), we have
  \bel{1209-b1}\ba{ll}&\ds
  A^2 z J_2(z)+2A\D zJ_2(z)\\
  \ns&\ds=\n\cdot V_{5}-2\n A^2\cdot\n\ell |\n z|^2-2A F |\n z|^2+2\n\cdot\(A^3\n\ell\) z^2+A^2Fz^2\\
 \ns&\ds\q+4\n\cdot(A\n\ell)|\D z|^2-4A\D\ell|\D z|^2+8A^2\n^2 \ell\n z\n z+4(\n A^2\cdot\n z)(\n\ell\cdot\n z)\\
 \ns&\ds\q-4\n\cdot\( A^2\n\ell\)|\n z|^2-8A\n^2\ell:\n^2 z\D z+R_2,
  \ea\ee

  where $V_{5}$ is given by (\ref{m0v-a0}) and
   \bel{1210-R2}\ba{ll}\ds
   R_2=-4 z\(\n z\cdot\n \ell\) \D A^2-4 z\(\n^2 \ell\n z\n A^2\)+2\n(A^2\D\ell)\cdot\n z z-2\n\(A F\) \cdot\n z z\\
   \ea\ee

 \ms

{\bf  Step 3. Estimation of $4(\n^2z\n\ell\n\ell)J_2(z)+\vec{E}\cdot\n zJ_2(z)$. }

\ms

By (\ref{0624-0a}) and (\ref{fg}), we have
   \bel{i14i22}\ba{ll}&\ds
4(\n^2 z\n\ell\n\ell)J_2(z)
\\
\ns&\ds=-16\n\cdot\[\(\n^2 z\n\ell\n\ell\)\n^2 z\n\ell\]+32\n^2\ell(\n^2 z\n\ell)(\n^2 z\n\ell)\\
\ns&\ds\q+8\n\cdot\(| \n^2 z\n\ell |^2\n\ell\)-8\D\ell | \n^2 z\n\ell |^2-16\n^2 z(\n^2\ell\n\ell)(\n^2 z\n\ell)\\
\ns&\ds\q-8\D\ell\D z\n^2 z\n\ell\n\ell-8\n\cdot\(A\n\ell|\n\ell\cdot\n z|^2\)+8\n\cdot\(A\n\ell\)|\n\ell\cdot\n z|^2\\
\ns&\ds\q+16A(\n^2\ell\n\ell\n z)(\n\ell\cdot\n z)+4\n\cdot\[F z(\n z\cdot\n\ell)\n\ell\]-4\(\n F\cdot\n\ell\)(\n z\cdot\n\ell) z\\
\ns&\ds\q-4 F\D\ell(\n z\cdot\n\ell) z-4 F(\n^2\ell\n\ell\n z) z-4 F|\n z\cdot\n\ell|^2.
\ea\ee
On the one hand, by (\ref{na}), we know that
 \bel{1209-b3}\left\{\ba{ll}\ds
-16\n^2 z\n^2\ell\n\ell(\n^2 z\n\ell)=-8\(\n^2 z\n A\)(\n^2 z\n\ell),\\
\ns\ds 16A(\n^2\ell\n\ell\n z)(\n\ell\cdot\n z)=8A(\n A\cdot\n z)(\n\ell\cdot\n z).
\ea\right.\ee
On the other hand,
\bel{yrr}\ba{ll}\ds
-8\D\ell\D z\n^2 z\n\ell\n\ell\\
\ns\ds=-8\n\cdot\[\D\ell\D z\(\n z\cdot\n\ell\)\n\ell\]+8\(\n\D\ell\cdot\n\ell\)\(\n z\cdot\n\ell\)\D z\\
\ns\ds\q+8|\D\ell|^2\(\n z\cdot\n\ell\)\D z+8\D\ell\D z\(\n^2\ell\n z\n\ell\)\\
\ns\ds\q+8\n\cdot\[\D\ell\(\n^2 z\cdot\n\ell\)\(\n z\cdot\n\ell\)\]-8\(\n^2 z\n\D\ell\n\ell\)\(\n\ell\cdot\n z\)\\
\ns\ds\q-8\D\ell\(\n^2 z:\n^2\ell\)\(\n\ell\cdot\n z\)-8\D\ell\(\n^2 z\cdot\n\ell\)\(\n^2\ell\cdot\n z\)-8\D\ell|\n^2 z\n\ell|^2.
\ea\ee
Recalling
  \bel{E-1}
  \vec{E}\cdot\n z=4(\n A+\D\ell\n\ell)\cdot\n z,
   \ee
then we have
 \bel{i15i22}\ba{ll}&\ds
\vec{E}\cdot\n zJ_2(z)\\
\ns&\ds=(\vec{E}\cdot\n z)\(-4 \n\ell\cdot\n\D z-4\n^2\ell:\n^2 z-2\D\ell\D z-4A\n\ell\cdot\n z+Fz\)\\
\ns&\ds=-4\n\cdot\((\vec{E}\cdot\n z)(\n^2 z\n\ell)\)+ 4\n^2 z\n\ell\n (\vec{E}\cdot\n z)-8\D \ell\D z(\n A\cdot\n z)\\
\ns&\ds\qq\q-8(\D\ell)^2\D z\n\ell\cdot\n z-8(\n A^2\cdot\n z)(\n\ell\cdot\n z)\\
\ns&\ds\qq\q-16A\D\ell|\n\ell\cdot\n z|^2+F\vec{E}\cdot \n z z.
 \ea\ee
Moreover, by (\ref{E-1}), it is easy to see that
 \bel{li-1}\ba{ll}\ds
 4\n^2 z\n\ell\n (\vec{E}\cdot\n z) \\
\ns\ds=16(\n^2 z\n\ell)(\n^2 z\n A)+16\D\ell |\n^2 z\n\ell|^2+16(\n^2 z\n\ell)(\n^2 A\n z)\\
 \ns\ds\q+16\D\ell(\n^2 z\n\ell)\cdot (\n^2\ell \n z)+16(\n^2 z\n\ell\n\D\ell)(\n\ell\cdot\n z).
  \ea\ee
Combining (\ref{i14i22})--(\ref{yrr}), (\ref{i15i22}) and (\ref{li-1}), we have
\bel{i1i22}\ba{ll}&\ds
4(\n^2z\n\ell\n\ell)J_2(z)+\vec{E}\cdot\n zJ_2(z)\\
\ns&\ds
=\n\cdot V_{6}+32\n^2\ell(\n^2 z\n\ell)(\n^2 z\n\ell)+8\[\n\cdot(A\n\ell)-2A\D\ell\]|\n\ell\cdot\n z|^2\\
\ns&\ds\q-8A(\n A\cdot\n z)(\n\ell\cdot\n z)+8(\n^2 z\n\ell)(\n^2 z\n A)-4 F|\n z\cdot\n\ell|^2+R_3,
\ea\ee
where
 \bel{1210-R3}\ba{ll}\ds
 R_3=8\(\n\D\ell\cdot\n\ell\)\(\n z\cdot\n\ell\)\D z-4\D\ell\D z\n A\cdot\n z+8(\n^2 z\n\ell\n\D\ell)(\n\ell\cdot\n z)\\
 \ns\ds\qq\q-8\D\ell\(\n^2 z:\n^2\ell\)\(\n\ell\cdot\n z\)+16(\n^2 z\n\ell)\cdot(\n^2 A\n z)+8\D\ell(\n^2 z\n\ell)\cdot (\n^2\ell \n z)\\
 \ns\ds\qq\q-4\(\n F\cdot\n\ell\)(\n z\cdot\n\ell) z-4 F\D\ell(\n z\cdot\n\ell) z-4 F(\n^2\ell\n\ell\n z) z+F\vec{E}\cdot \n z z.
 \ea\ee

 \ms

 {\bf Step 4: Estimation of ``$\ds J_1(z)\cP_3(z)$". } First,
   \bel{i11i23-0}\ba{ll}&\ds
\D^2 z\cP_3(z)\ds=\D^2 z \(\Phi z+ Z\n^2 z:\n^2\eta\)\\
\ns&\ds= \n\cdot\(\Phi \n\D z z-\Phi\D z\n z\)+\Phi |\D z |^2-\n \Phi \cdot\n\D z z+\n \Phi\cdot\n z \D z\\
\ns&\ds\q+\n\cdot\(( Z\n^2 z:\n^2\eta)\n\D z\)-\sum_{i, j, k=1}^n\(Z\eta_{x_ix_k}\)_{x_j}\D z_{x_j}z_{x_ix_k}\\
\ns&\ds\q-\sum_{i, j, k=1}^n\( Z\eta_{x_ix_k}\D z_{x_j}z_{x_ix_j}\)_{x_k}+\sum_{i, j, k=1}^n\( Z\eta_{x_ix_k}\)_{x_k}\D z_{x_j}z_{x_ix_j}\\
\ns&\ds\q+\sum_{i, j, k=1}^n\( Z\eta_{x_ix_k}\D z_{x_k}z_{x_ix_j}\)_{x_j}-\sum_{i, j, k=1}^n\( Z\eta_{x_ix_k}\)_{x_j}\D z_{x_k}z_{x_ix_j}-Z\n^2\eta\n\D z\n\D z.
\ea\ee

 Next,
   \bel{i12i24}\ba{ll}
A^2 zJ_3(z)&\ds=A^2 z \(\Phi z+ Z\n^2 z:\n^2\eta\)\\
\ns&\ds=A^2\Phi  z^2+\n\cdot\(A^2{ Z}(\n^2\eta\n z) z\)\\
\ns&\ds\q-\(\n^2\eta\n(A^2{ Z})\n z \)z-A^2{ Z}\n^2\eta\n z\n z-A^2{ Z}(\n\D\eta\cdot\n z) z
\ea\ee
and
 \bel{i14i24}\ba{ll}&\ds
2A\D zJ_3(z)\\
\ns&\ds= 2A\D z \(\Phi z+ Z\n^2 z:\n^2\eta\)\\
 \ns&\ds=\n\cdot\[2A\Phi z\n z\]-2\n\(A\Phi\)\cdot \n z z-2A\Phi |\n z|^2 +2A Z(\n^2 z:\n^2\eta)\D z
 \ea\ee
 Further, by (\ref{E-1}), we have
  \bel{i1-4}\ba{ll}&\ds
4\((\n^2 z\n\ell)\cdot\n\ell\) J_3(z)\\
\ns&\ds=4\Phi z(\n^2 z\n\ell\n\ell)+4 Z(\n^2 z:\n^2\eta)(\n^2 z\n\ell\n\ell)\\
\ns&\ds=4\n\cdot\[\Phi z(\n z\cdot\n\ell)\n\ell\]-4\(\n\Phi\cdot\n\ell\)(\n z\cdot\n\ell) z-4\Phi\D\ell(\n z\cdot\n\ell) z\\
\ns&\ds\q-4\Phi(\n^2\ell\n\ell\n z) z-4\Phi|\n z\cdot\n\ell|^2+4 Z(\n^2 z:\n^2\eta)(\n^2 z\n\ell\n\ell)
 \ea\ee
 and
 \bel{i15i24}\ba{ll}
\vec{E}\cdot\n z J_3(z)&\ds=\Phi\vec{E}\cdot\n z z+ Z\(\n^2 z:\n^2\eta\)(\vec{E}\cdot\n z).
 \ea\ee

Combining (\ref{i11i23-0})--(\ref{i15i24}), we have
 \bel{1210-a0}\ba{ll}&\ds
 J_1(z)J_3(z)\\
 \ns&\ds=\n\cdot V_{7}+\Phi |\D z |^2-2A\Phi |\n z|^2-A^2{ Z}(\n^2\eta\n z)\n z+A^2\Phi  z^2-Z\n^2\eta\n\D z\n\D z\\
\ns&\ds\q+2A Z(\n^2 z:\n^2\eta)\D z-4\Phi|\n z\cdot\n\ell|^2 +4 Z(\n^2 z:\n^2\eta)(\n^2 z\n\ell\n\ell)+\mathcal{R}
 \ea\ee
where $\mathcal{R}$ is given by \eqref{1210-R4}.

 \ms

 {\bf Step 5: End of the proof. }In order to finish the proof, it remains to estimate the remainder terms $R_i,i=1,2,3$. In view of  condition \eqref{cd1}, combing \eqref{0421-a} with \eqref{1210-R1},  we get
\bel{R1-es}\ba{ll}\ds
|R_1|\leq C\[\lambda\mu^3\varphi|\n \D z\cdot \n \eta||\D z|+\lambda\mu^2\varphi|\n \D z||\D z|+\lambda\mu^4\varphi|\n^2 z||\D z| \\
\ns\ds\qq\q+\lambda\mu^3\varphi|\n \D z\cdot \n \eta||\n^2 z|+\lambda\mu^2\varphi|\n \D z||\n^2 z|+\lambda^3\mu^5\varphi^3|\n\eta\cdot\n z||\D z|\\
\ns\ds\qq\q+\lambda^3\mu^4\varphi^3|\n z||\D z|+\lambda^3\mu^6\varphi^3|z||\n^2 z|+ \lambda^3\mu^5\varphi^3|\n \D z\cdot \n \eta||z| \\
\ns\ds\qq\q +\lambda^3\mu^5\varphi^3|\n\eta\cdot\n z||\n^2 z|+\lambda^3\mu^4\varphi^3|\n z||\n^2 z| +\lambda^3\mu^4\varphi^3|\n \D z|| z|
\]
\ea\ee

Recalling the definition \eqref{A11-12-13} of $A_i,i=1,2,3,4$, using inequality $ab\geq -\frac{1}{2}(a^2+b^2)$, we have
\bel{R1-es1}\ba{ll}\ds
R_1\geq -C \(A_1+A_2+A_4\).
\ea\ee
Similarly, we have
\bel{1210-R21} \ba{ll}\ds
|R_2|\leq C\[ \lambda^3\mu^5\varphi^3 |z||\n\eta\cdot\n z|+\lambda^3\mu^4\varphi^3|z||\n z|
\],
\ea\ee
and
 \bel{1210-R031}\ba{ll}\ds
 |R_3|\leq C\[\lambda^3\mu^5\varphi^3|\n z\cdot \n \eta||\D z|+\lambda^3\mu^4\varphi^3|\n z||\D z|+\lambda^3\mu^5\varphi^3|\n z\cdot \n \eta||\n^2 z| \\
  \ns\ds\qq\q+\lambda^3\mu^4\varphi^3|\n z||\n^2 z| +\lambda^4\mu^6\varphi^4|\n z\cdot \n \eta|| z|+\lambda^4\mu^5\varphi^4|\n z|| z|
 \],
\ea\ee
 which give
 \bel{1210-R431}\ba{ll}\ds
R_2\geq -C A_1, \q R_3\geq -C(A_1+A_2).
\ea\ee

 Hence together (\ref{R1-es1}) with (\ref{1210-R431}), we have
    \bel{1210-b0}\ba{ll}\ds
  R_1+R_2+R_3\ge -C(A_1+A_2+A_3).
 \ea\ee

Finally, combining (\ref{1209-06}), (\ref{1209-b1}), (\ref{i1i22}), (\ref{1210-a0}) and \eqref{1210-b0}, we get the desired result immediately. \endpf


\begin{thebibliography}{99}

\bibitem{AdamsBook}
Robert~A. Adams and John J.~F. Fournier.
\newblock {\em Sobolev spaces}, volume 140 of {\em Pure and Applied Mathematics
  (Amsterdam)}.
\newblock Elsevier/Academic Press, Amsterdam, second edition, 2003.

\bibitem{adn}
Shmuel Agmon, Avron Douglis, and Louis Nirenberg.
\newblock Estimates near the boundary for solutions of elliptic partial
  differential equations satisfying general boundary conditions. {I}.
\newblock {\em Commun. Pure Appl. Math.}, 12:623--727, 1959.

\bibitem{BL15}
Mourad Bellassoued and J{\'e}r{\^o}me Le~Rousseau.
\newblock Carleman estimates for elliptic operators with complex coefficients.
  {II}: {Transmission} problems.
\newblock {\em J. Math. Pures Appl. (9)}, 115:127--186, 2018.

\bibitem{BZ}
N.~Burq and C.~Zuily.
\newblock Concentration of {Laplace} eigenfunctions and stabilization of weakly
  damped wave equation.
\newblock {\em Commun. Math. Phys.}, 345(3):1055--1076, 2016.

\bibitem{burq98}
Nicolas Burq.
\newblock D\'ecroissance de l'\'energie locale de l'\'equation des ondes pour
  le probl\`eme ext\'erieur et absence de r\'esonance au voisinage du r\'eel.
\newblock {\em Acta Math.}, 180(1):1--29, 1998.

\bibitem{c}
T.~Carleman.
\newblock Sur une probl\`eme d'unicit\'e pour les syst\`emes d'\'equations aux
  d\'eriv\'ees partielles \`a deux variables ind\'ependantes.
\newblock {\em Ark. Mat. Astr. Fys.}, 26B(17):1--9, 1939.

\bibitem{CM}
Eduardo Cerpa and Alberto Mercado.
\newblock Local exact controllability to the trajectories of the 1-{D}
  {Kuramoto}-{Sivashinsky} equation.
\newblock {\em J. Differ. Equations}, 250(4):2024--2044, 2011.

\bibitem{coron}
Jean-Michel Coron.
\newblock {\em {Control and nonlinearity}}.
\newblock Amer Mathematical Society, 2007.

\bibitem{D}
Thomas Duyckaerts.
\newblock Optimal decay rates of the energy of a hyperbolic-parabolic system
  coupled by an interface.
\newblock {\em Asymptotic Anal.}, 51(1):17--45, 2007.

\bibitem{FLZ}
Xiaoyu Fu, Qi~L{\"u}, and Xu~Zhang.
\newblock {\em Carleman estimates for second order partial differential
  operators and applications. {A} unified approach}.
\newblock SpringerBriefs Math. Cham: Springer, 2019.

\bibitem{FI}
A.~V. Fursikov and O.~Yu. Imanuvilov.
\newblock {\em Controllability of evolution equations}, volume~34 of {\em Lect.
  Notes Ser., Seoul}.
\newblock Seoul: Seoul National Univ., 1996.

\bibitem{GP0}
Peng Gao.
\newblock Insensitizing controls for the {Cahn}-{Hilliard} type equation.
\newblock {\em Electron. J. Qual. Theory Differ. Equ.}, 2014:22, 2014.


\bibitem{GP}
Peng Gao.
\newblock A new global {Carleman} estimate for the one-dimensional
  {Kuramoto}-{Sivashinsky} equation and applications to exact controllability
  to the trajectories and an inverse problem.
\newblock {\em Nonlinear Anal., Theory Methods Appl., Ser. A, Theory Methods},
  117:133--147, 2015.

\bibitem{ggs}
Filippo Gazzola, Hans-Christoph Grunau, and Guido Sweers.
\newblock {\em Polyharmonic boundary value problems. {Positivity} preserving
  and nonlinear higher order elliptic equations in bounded domains}, volume
  1991 of {\em Lect. Notes Math.}
\newblock Berlin: Springer, 2010.

\bibitem{GK}
Shannon~Guerrero and Karim~Kassab.
\newblock Carleman estimate and null controllability of a fourth order
  parabolic equation in dimension {{\($N\geq 2$\)}}.
\newblock {\em J. Math. Pures Appl. (9)}, 121:135--161, 2019.

\bibitem{Hormander63}
Lars H{\"o}rmander.
\newblock {\em Linear partial differential operators}.
\newblock Die Grundlehren der mathematischen Wissenschaften, Bd. 116. Academic
  Press, Inc., Publishers, New York; Springer-Verlag,
  Berlin-G\"ottingen-Heidelberg, 1963.

\bibitem{Hormander}
Lars H{\"o}rmander.
\newblock Uniqueness theorems for second order elliptic differential equations.
\newblock {\em Commun. Partial Differ. Equations}, 8:21--64, 1983.

\bibitem{HK}
Xinchi Huang and Atsushi Kawamoto.
\newblock Inverse problems for a half-order time-fractional diffusion equation
  in arbitrary dimension by {Carleman} estimates.
\newblock {\em Inverse Probl. Imaging}, 16(1):39--67, 2022.

\bibitem{Isakov2}
Victor Isakov.
\newblock {\em Inverse problems for partial differential equations}, volume 127
  of {\em Appl. Math. Sci.}
\newblock New York, NY: Springer, 2nd ed. edition, 2006.

\bibitem{Kassab}
Karim Kassab.
\newblock Null controllability of semi-linear fourth order parabolic equations.
\newblock {\em J. Math. Pures Appl. (9)}, 136:279--312, 2020.

\bibitem{Kn}
Carlos~E. Kenig.
\newblock Carleman estimates, uniform {Sobolev} inequalities for second-order
  differential operators, and unique continuation theorems.
\newblock Proc. {Int}. {Congr}. {Math}., {Berkeley}/{Calif}. 1986, {Vol}. 2,
  948-960 (1987)., 1987.

\bibitem{Kli}
Michael~V. Klibanov.
\newblock Carleman estimates for global uniqueness, stability and numerical
  methods for coefficient inverse problems.
\newblock {\em J. Inverse Ill-Posed Probl.}, 21(4):477--560, 2013.

\bibitem{JL20}
J{\'e}r{\^o}me Le~Rousseau and Luc Robbiano.
\newblock Spectral inequality and resolvent estimate for the bi-{Laplace}
  operator.
\newblock {\em J. Eur. Math. Soc. (JEMS)}, 22(4):1003--1094, 2020.

\bibitem{Lionscontrolbook}
Jacques-Louis Lions.
\newblock {\em {Contr{\^o}labilit{\'e} Exacte, Stabilization et Perturbations
  de Syst\`emes Distribu{\'e}es, Tom 1}}.
\newblock Masson, RMA, 1988.

\bibitem{XCY}
Xiang Xu, Jin Cheng, and Masahiro Yamamoto.
\newblock Carleman estimate for a fractional diffusion equation with half order
  and application.
\newblock {\em Appl. Anal.}, 90(9-10):1355--1371, 2011.

\bibitem{Y}
Masahiro Yamamoto.
\newblock Carleman estimates for parabolic equations and applications.
\newblock {\em Inverse Probl.}, 25(12):75, 2009.


\bibitem{Zheng}
Chuang Zheng.
\newblock Inverse problems for the fourth order {Schr{\"o}dinger} equation on a
  finite domain.
\newblock {\em Math. Control Relat. Fields}, 5(1):177--189, 2015.

\bibitem{ZZ}
Zhongcheng Zhou.
\newblock Observability estimate and null controllability for one-dimensional
  fourth order parabolic equation.
\newblock {\em Taiwanese J. Math.}, 16(6):1991--2017, 2012.

\bibitem{Zuily83}
Claude Zuily.
\newblock {\em Uniqueness and non-uniqueness in the {Cauchy} problem},
  volume~33 of {\em Prog. Math.}
\newblock Birkh{\"a}user, Cham, 1983.

     \end{thebibliography}
\end{document}